\documentclass[11pt]{article}

\usepackage{latexsym}
\usepackage{amssymb}
\usepackage{graphicx}

\newtheorem{Theorem}{Theorem}[part]

\newtheorem{Remark}{Remark}[part]

\def \Prod{\displaystyle\prod}

\def \Frac{\displaystyle\frac}

\def \N{\mathbb{N}}
\def \R{\mathbb{R}}

\def \E{\mathbb{E}}
\def \F{\mathbb{F}}

\def \P{\mathbb{P}}
\def \Q{\mathbb{Q}}

\def \Ac{{\cal A}}
\def \Bc{{\cal B}}

\def \Dc{{\cal D}}
\def \Ec{{\cal E}}
\def \Fc{{\cal F}}

\def \Pc{{\cal P}}

\def \Nc{{\cal N}}

\def \Xc{{\cal X}}

\def \eps{\varepsilon}

\def \ep{\hbox{ }\hfill$\Box$}

\def\Dt#1{\Frac{\partial #1}{\partial t}}

\def\reff#1{{\rm(\ref{#1})}}

\def\beqs{\begin{eqnarray*}}
\def\enqs{\end{eqnarray*}}
\def\beq{\begin{eqnarray}}
\def\enq{\end{eqnarray}}

\begin{document}

\title{Some applications and methods of large deviations \\ in finance and insurance
\thanks{Lectures from a Bachelier course at Institut Henri Poincar\'e, Paris, February 2007. 
I would like to thank the participants for their interest and comments.} }

\author{Huy{\^e}n PHAM
             \\\small  Laboratoire de Probabilit{\'e}s et
             \\\small  Mod{\`e}les Al{\'e}atoires
             \\\small  CNRS, UMR 7599
             \\\small  Universit{\'e} Paris 7
             \\\small  e-mail: pham@math.jussieu.fr
             \\\small  and Institut Universitaire de France}

\maketitle

\begin{abstract}
In these notes, we present  some methods and applications of large deviations  to finance and insurance.  
We begin with  the classical ruin problem related to the Cramer's theorem and give en extension to  an insurance model with investment in stock market.  We then describe  how large deviation approximation  and importance sampling are used  in rare event simulation for option pricing. 
We finally focus on large deviations methods in risk management for the 
estimation of  large portfolio losses  in credit risk and portfolio performance in market investment.  
\end{abstract}

\vspace{7mm}

\noindent {\bf Key words~:}   large deviations,  ruin problem, importance sampling, rare event simulation, exit probability, credit risk, portfolio performance.

\vspace{5mm}

\noindent {\bf MSC Classification (2000)~:}  60F10, 62P05,  65C05,  91B28,  91B30.

\newpage

\tableofcontents

\newpage

\section{Introduction}

\setcounter{equation}{0} \setcounter{Assumption}{0}
\setcounter{Theorem}{0} \setcounter{Proposition}{0}
\setcounter{Corollary}{0} \setcounter{Lemma}{0}
\setcounter{Definition}{0} \setcounter{Remark}{0}

The area of large deviations is a set of asymptotic results on rare event probabilities and a set of methods to derive such results.  Large deviations is a very active area in applied probability, and questions related to extremal events both in finance and insurance applications, play an increasingly important role. 
For instance, recent applications of interest  concern  ruin probabilities in risk theory,  risk management, 
option pricing, and simulation of quantities involved in this context.

Large deviations appear historically in insurance mathematics with the ruin probability estimation problem within the classical 
Cramer-Lundberg model. The problem was then subsequently extended to more general  models involving for example L\'evy 
risk processes.  In finance, large deviations  arise in various contexts.   They occur  in risk management for the computation 
of  the probability of large losses  of a portfolio subject to market risk or the default probabilities of a portfolio under credit risk.  Large deviations 
methods are largely used   in rare events simulation and so appear  naturally in the approximation  of option pricing, in particular for barrier options and far from the money options.

We illustrate our purpose with the following toy example.  Let $X$ be a (real-valued) random variable, and consider the problem of computing or estimating 
$\P[X > \ell]$, the probability that $X$ exceeds some level $\ell$.  
In finance, we may  think of $X$ as the loss of a portfolio subject to  market or credit risk, and 
we are interested in the probability of large loss or default probability.  The r.v. $X$ may also correspond to the terminal value of a stock price, and the quantity $\P[X > \ell]$ appears typically in the computation of a call or barrier option, with a small probability of payoff  when the option is far from the money or the barrier $\ell$ is large.  To estimate $p$ $=$ $\P[X > \ell]$, a basic technique is Monte Carlo simulation~: generate $n$ independent copies  
$X_1,\ldots,X_n$ of $X$, and use the sample mean~:
\beqs
\bar S_n &=& \frac{1}{n} \sum_{i=1}^n Y_i, \;\;\;\; \mbox{ with  } \; Y_i = 1_{X_i > \ell}. 
\enqs
The convergence of this estimate (when $n$ $\rightarrow$ $\infty$)  follows from the law of large numbers, while the standard rate of convergence is 
given, via the central limit theorem, in terms of the variance $v$ $=$ $p(1-p)$ of $Y_i$~: 
\beqs
\P\big[ | \bar S_n - p | \geq \frac{a}{\sqrt{n}} \big] & \rightarrow & 2 \Phi\big(- \frac{a}{\sqrt{v}} \big),
\enqs
where $\Phi$ is the cumulative distribution function of the standard normal law.  Furthermore, the convergence of the estimator $\bar S_n$ is  precised with the large deviation result, known here as the Cramer's theorem, which is concerned with approximation of rare event probabilities $\P[ \bar S_n \in A]$, and typically states that
\beqs
\P \big[ | \bar S_n - p | \geq  a \big] & \simeq & C e^{-\gamma n},
\enqs
for some constants $C$ and $\gamma$. 

Let us now turn again to the estimation of $p$ $=$ $\P[X>\ell]$. As mentioned above, the rate of convergence of the naive estimator $\bar S_n$ is 
determined by~: 
\beqs
{\rm Var}(\bar S_n) &=&  \frac{{\rm Var}(1_{X>\ell})}{n} \; = \;  \frac{p(1- p)}{n},
\enqs
and the relative error is 
\beqs
\mbox{relative error} &=& \frac{\mbox{standard deviation of } \;  \bar S_n}{\mbox{mean of } \bar S_n} \; = \; \frac{ \sqrt{p(1- p)}}{\sqrt{n} p}. 
\enqs 
Hence, if $p$ $=$ $\P[X >\ell]$ is small, and since $\sqrt{p-p^2}/p$ $\rightarrow$ $\infty$ as $p$ goes to zero, we see that a large sample size (i.e. $n$) is required for the estimator to achieve a  reasonable relative error bound.  This is a common occurence when estimating rare events.  In order to improve 
the estimate of  the tail probability $\P[X >\ell]$,  one is tempted to use importance sampling to reduce variance, and hence  speed up the computation by requiring fewer samples.  This consists basically in changing measures  to try to give more weight to ``important" outcomes,  
(increase the default probability).  Since large deviations theory also deals with rare events, we can see  its strong link with importance sampling.  
 
 To make the idea concrete, consider again the problem of estimating $p$ $=$ $\P[X>\ell]$, and suppose that $X$ has distribution $\mu(dx)$. Let us look at an alternative sampling distribution $\nu(dx)$  absolutely continuous with respect to $\mu$, with density $f(x)$ $=$ $d\nu/d\mu(x)$.  The 
tail probability can then be written as~:
\beqs
p \; = \;  \P[ X > \ell ] &=&  \int 1_{x > \ell} \mu(dx)  \; = \; \int  1_{x > \ell} \phi(x) \nu(dx) \; = \; \E_\nu [1_{X > \ell} \phi(X) ],
\enqs 
where $\phi$ $=$ $1/f$, and $\E_\nu$ denotes  the expectation under the measure $\nu$. By generating i.i.d. samples  $\tilde X_1, \ldots,\tilde X_n,\ldots$ with distribution $\nu$, we have then an alternative unbiased and convergent estimate of $p$ with 
\beqs
\tilde S_n &=&   \frac{1}{n} \sum_{i=1}^n 1_{\tilde X_i > \ell} \phi(\tilde X_i),
\enqs
and  whose  rate of convergence is determined by
\beqs
{\rm Var}_\nu(\bar S_n) &=& \frac{1}{n} \int \big(1_{x > \ell}  - p f(x) \big)^2 \phi^2(x) \nu(dx).   
\enqs 
The minimization of this quantity over all possible $\nu$ (or $f$) leads to a  zero variance with the choice of a density $f(x)$ $=$ $1_{x>\ell}/p$. This is of course only a theoretical result since it requires the knowledge of $p$, the very thing we want to estimate! However,  by noting that in this case 
$\nu(dx)$ $=$ $f(x)\mu(dx)$ $=$ $1_{x>\ell} \mu(dx)/\P[X > \ell]$ is nothing else than the conditional distribution of $X$ given $\{X>\ell\}$, this suggests to use an importance sampling change of measure that makes the rare event  $\{X > \ell\}$ more likely.  This method of suitable change of measure is 
also the key  step in proving large deviation results.

The plan of these lectures is the following. In Section 2, we give some basic tools and results on large deviations, in particular the  most classical result on large deviations,  Cramer's theorem. We illustrate in Section 3 the first applications of large deviations to ruin problems in insurance industry, and give some extension to an insurance model with financial investment opportunity.  
Section 4  is concerned with the large deviations approximation for rare event simulation in option pricing, and we shall use 
asymptotic results from large deviations theory~: Fredilin-Wentzell theory on sample path large deviations, and Varadhan's integral principle. Finally, Section 5 is devoted to applications of large deviations in risk management,  where we use conditional and control variants of the Ellis-Gartner theorem.

\section{Basic tools and results  on large deviations}

\setcounter{equation}{0} \setcounter{Assumption}{0}
\setcounter{Theorem}{0} \setcounter{Proposition}{0}
\setcounter{Corollary}{0} \setcounter{Lemma}{0}
\setcounter{Definition}{0} \setcounter{Remark}{0}

\subsection{Laplace function and exponential change of measures}

If $X$ is a (real-valued) random variable on $(\Omega,\Fc)$ with probability distribution $\mu(dx)$, the cumulant generating function  (c.g.f.) 
of $\mu$ is the logarithm of the 
Laplace function of $X$, i.e.~:
\beqs
\Gamma(\theta) &=& \ln \E[e^{\theta X}] \; = \; \ln \int e^{\theta x} \mu(dx)  \; \in \; (-\infty,\infty], \;\;\;\; \theta \in \R. 
\enqs  
Notice that  $\Gamma(0)$ $=$ $0$, and $\Gamma$ is  convex by 
H\"older inequality.  We denote ${\cal D}(\Gamma)$ $=$ $\{\theta \in \R~: \Gamma(\theta) < \infty\}$, and for any 
$\theta$ $\in$ $\Dc(\Gamma)$, we define a probability measure $\mu_\theta$ on $\R$ by~: 
\beq \label{mutheta}
\mu_\theta(dx)   &=&   \exp(\theta x - \Gamma(\theta)) \mu(dx).  
\enq
Suppose  that $X_1,\ldots,X_n, \ldots,$  is an i.i.d.  sequence of  
random variables  with distribution $\mu$ and consider the new probability measure  $\P_\theta$ on 
$(\Omega,\Fc)$ with likelihood ratio evaluated at $(X_1,\ldots,X_n)$, $n$ $\in$ $\N^*$,  by~: 
\beq \label{ratioptheta}
\frac{d\P_\theta}{d\P}(X_1,\ldots,X_n) &=& \Prod_{i=1}^n \frac{d\mu_\theta}{d\mu}(X_i) \; = \; \exp\Big(\theta\sum_{i=1}^n X_i - n \Gamma(\theta)\Big). 
\enq
By denoting $\E_\theta$ the corresponding expectation under $\P_\theta$,  formula \reff{ratioptheta} means that for all $n$ $\in$ $\N^*$, 
\beq \label{changeesp}
\E\Big[ f(X_1,\ldots,X_n) \Big] &=& \E_\theta\Big[  f(X_1,\ldots,X_n)  \exp\Big(-\theta\sum_{i=1}^n X_i + n \Gamma(\theta)\Big)  \Big],
\enq
for all  Borel functions $f$ for which the expectation on the l.h.s. of \reff{changeesp} is finite.  Moreover,  the random variables $X_1,\ldots,X_n$, $n$ $\in$ 
$\N^*$,  are i.i.d.  with probability distribution $\mu_\theta$ under $\P_\theta$.   Actually, the relation \reff{changeesp} extends from a fixed number of 
steps $n$ to a random number of steps, provided the random horizon is a stopping time.  More precisely, if $\tau$ is a stopping time in $\N$  for 
$X_1,\ldots,X_n,\ldots$, i.e. the event $\{\tau  < n \}$ is measurable with respect to the algebra generated by $\{X_1,\ldots,X_n\}$ for all $n$,  then 
\beq \label{changeesptau}
\E\Big[ f(X_1,\ldots,X_\tau) 1_{\tau < \infty}  \Big] &=& \E_\theta\Big[  f(X_1,\ldots,X_\tau)  
\exp\Big(-\theta\sum_{i=1}^\tau X_i + \tau \Gamma(\theta)\Big)  1_{\tau < \infty} \Big], 
\enq
for all  Borel functions $f$ for which the expectation on the l.h.s. of \reff{changeesptau} is finite.

The cumulant generating function $\Gamma$ records some useful information on the probability distributions $\mu_\theta$.  For example, 
$\Gamma'(\theta)$ is the mean of $\mu_\theta$. Indeed, for any $\theta$ in the interior of $\Dc(\Gamma)$, differentiation yields by dominated convergence~: 
\beq \label{changemean}
\Gamma'(\theta) &=& \frac{\E[Xe^{\theta X}]}{\E[e^{\theta X}]} \; = \; \E\big[X \exp\big(\theta X - \Gamma(\theta) \big)\big] \; = \; \E_\theta[X]. 
\enq
A similar calculation shows that $\Gamma''(\theta)$ is the variance of $\mu_\theta$. 
Notice in particular that if $0$ lies in the interior of $\Dc(\Gamma)$, then $\Gamma'(0)$ $=$ $\E[X]$ and $\Gamma''(0)$ $=$ $Var(X)$. 

\vspace{2mm}

\noindent {\bf Bernoulli distribution}

\noindent Let $\mu$ the Bernoulli distribution of parameter $p$. Its c.g.f. is given by
\beqs
\Gamma(\theta) &=&  \ln (1-p + p e^\theta). 
\enqs 
A direct simple algebra calculation shows that  $\mu_\theta$  is the Bernoulli distribution of  parameter 
$p e^{\theta}/(1-p+pe^\theta)$.

\vspace{2mm}

\noindent {\bf Poisson distribution}

\noindent Let $\mu$ the Poisson distribution of intensity $\lambda$. Its c.g.f. is given by
\beqs
\Gamma(\theta) &=&  \lambda(e^{\theta} - 1). 
\enqs 
A direct simple algebra calculation shows that  $\mu_\theta$  is the Poisson distribution of intensity $\lambda e^{\theta}$. 
Hence, the effect of the change of probability measure $\P_\theta$ is to multiply  the intensity by a factor  $e^\theta$.

\vspace{2mm}

\noindent {\bf Normal distribution}

\noindent Let $\mu$ the normal distribution $\Nc(0,\sigma^2)$, whose c.g.f. is given by~: 
\beqs
\Gamma(\theta) &=&  \frac{\theta^2\sigma^2}{2}. 
\enqs
A direct simple algebra calculation shows that $\mu_\theta$ is the normal distribution $\Nc(\theta\sigma^2,\sigma^2)$.  Hence, if $X_1,\ldots,X_n$ are 
i.i.d.  with normal distribution $\Nc(0,\sigma^2)$, then under the change of measure $\P_\theta$ with likelihood ratio~:
\beqs
\frac{d\P_\theta}{d\P}(X_1,\ldots,X_n) &=&   \exp\Big(\theta\sum_{i=1}^n X_i - n  \frac{\theta^2\sigma^2}{2} \Big),
\enqs
the random variables $X_1,\ldots,X_n$ are i.i.d. with  normal distribution $\Nc(\theta\sigma^2,\sigma^2)$~:  the effect of $\P_\theta$ is to change the mean 
of $X_i$ from $0$ to $\theta\sigma^2$.  This result  can be interpreted as  the finite-dimensional version of Girsanov's theorem.

\vspace{2mm}

\noindent {\bf Exponential distribution}

\noindent Let $\mu$ the exponential distribution of intensity $\lambda$. Its c.g.f. is given by
\beqs
\Gamma(\theta) &=& \left\{ \begin{array}{cc}
                                       \ln\big(\frac{\lambda}{\lambda - \theta}\big), & \;\; \theta < \lambda \\
                                       \infty, & \;\; \theta \geq \lambda 
                                       \end{array}
                                       \right.
\enqs 
A direct simple algebra calculation shows that for $\theta$ $<$ $\lambda$, $\mu_\theta$  is the exponential distribution of intensity $\lambda-\theta$. 
Hence, the effect of the change of probability measure $\P_\theta$ is to shift the intensity from $\lambda$ to $\lambda-\theta$.

\subsection{Cramer's theorem}

The most classical result in large deviations area is Cramer's theorem. This concerns large deviations associated with the empirical mean of i.i.d. random variables valued in a finite-dimensional space. We do not state the Cramer's theorem in whole generality.  Our purpose is to 
put emphasis on the methods used to derive such result.  For simplicity, we consider the case of real-valued  i.i.d. random variables  $X_i$ with (nondegenerate) probability distribution $\mu$ of  finite mean  $\E X_1$ $=$ $\int x\mu(dx)$ $<$ $\infty$, 
and we  introduce the random walk $S_n$ $=$ $\sum_{i=1}^n X_i$.  
It is well-known by the law of large numbers that  the empirical mean $S_n/n$ converges in probability to $\bar x$ $=$ $\E X_1$, i.e. $\lim_n \P[S_n/n \in (\bar x-\eps,\bar x+\eps)]$ $=$ $1$ for all $\eps$ $>$ $0$.  Notice also,  by the central limit theorem that $\lim_n \P[S_n/n \in [\bar x,\bar x+\eps)]$ $=$ $1/2$ for all $\eps$ $>$ $0$. 
Large deviations results focus on  asymptotics for probabilities of rare events, for example of  the form 
$\P\big[\frac{S_n}{n} \geq x \big]$ for $x$ $>$ $\E X_1$, and state that 
\beqs
\P\big[\frac{S_n}{n} \geq x \big] & \simeq &   C e^{-\gamma x}, 
\enqs  
for some  constants $C$ and  $\gamma$ to be precised later. 
The symbol $\simeq$ means that the ratio  is one in the limit (here when $n$ goes to infinity).  The rate of 
convergence is characterized by the Fenchel-Legendre transform of  the c.g.f. $\Gamma$ of $X_1$~:
\beqs
\Gamma^*(x) &=& \sup_{\theta \in \R}\big[ \theta x - \Gamma(\theta) \big] \; \; \in \; [0,\infty], \;\;\;  x \in \R. 
\enqs
As supremum  of affine functions, $\Gamma^*$ is convex. The sup in the definition  of $\Gamma^*$ can be evaluated by differentiation~: for $x$ $\in$ $\R$, 
if $\theta$ $=$ $\theta(x)$ is solution to the saddle-point equation, $x$ $=$ $\Gamma'(\theta)$,  then $\Gamma^*(x)$ $=$ $\theta x - \Gamma(\theta)$.  
Notice, from \reff{changemean}, that the exponential change of measure $\P_\theta$ put the expectation of $X_1$ to $x$.  Actually, exponential change of measure is a key tool in large deviations methods.  The idea is to select  a measure under which the rare event is no longer rare, so that the rate of 
decrease of the original probability is given by the rate of decrease of the likelihood ratio. This particular change of measure is intended to approximate 
the most likely way for the rare event to occur.  

By Jensen's inequality, we show  that  $\Gamma^*(\E X_1)$ $=$ $0$.  This implies that for all 
$x$ $\geq$ $\E X_1$, $\Gamma^*(x)$ $=$  $\sup_{\theta \geq 0}\big[ \theta x - \Gamma(\theta) \big]$, and so $\Gamma^*$ is nondecreasing on 
$[\E X_1,\infty)$.

\begin{Theorem} (Cramer's theorem) 

\noindent For any $x$ $\geq$ $\E X_1$, we have 
\beq \label{cramer}
\lim_{n\rightarrow\infty} \frac{1}{n} \ln \P \big[\frac{S_n}{n} \geq x \big]  &=& -  \Gamma^*(x) \; = \; - \inf_{y\geq x} \Gamma^*(y). 
\enq
\end{Theorem}
{\bf Proof.} 1) {\it Upper bound.}  The main step in the upper bound  $\leq$ of \reff{cramer} is based on Chebichev inequality combined with the i.i.d. assumption on the $X_i$~:  
\beqs
\P \big[\frac{S_n}{n} \geq x \big]  \; = \;  \E\big[ 1_{\frac{S_n}{n} \geq x }\big]  & \leq & \E\big[ e^{\theta(S_n - nx)}\big]  \; 
= \;  \exp\big( n\Gamma(\theta) - \theta nx\big), \;\;\; \forall \theta \geq 0. 
\enqs
By taking the infimum over $\theta$ $\geq$ $0$, and since 
$\Gamma^*(x)$ $=$ $\sup_{\theta\geq 0}[\theta x - \Gamma(\theta)]$ for $x$ $\geq$ $\E X_1$,  we then obtain
\beqs
\P\big[\frac{S_n}{n} \geq  x \big] & \leq & \exp\big(-n\Gamma^*(x)\big).  
\enqs
and so in particular the upper bound $\leq$  of  \reff{cramer}. 

\vspace{1mm}

\noindent 2) {\it Lower bound.}   Since $\P \big[\frac{S_n}{n} \geq x \big]$ $\geq$  $\P \big[\frac{S_n}{n} \in  [x,x+\eps) \big]$,  for all 
$\eps$ $>$ $0$,   it suffices  to show that 
\beq \label{interlowery}
\lim_{\eps\rightarrow 0} \liminf_{n\rightarrow\infty} \frac{1}{n} \ln \P \Big[\frac{S_n}{n} \in  [x,x+\eps) \Big]  & \geq & - \Gamma^*(x). 
\enq
Suppose  that $\mu$ is supported on a bounded support so that $\Gamma$ is finite everywhere.  Suppose  first that  
there exists a solution $\theta$ $=$ $\theta(x)$ $>$ $0$ to the saddle-point equation~:  $\Gamma'(\theta)$ $=$ $x$, i.e. attaining the supremum in 
$\Gamma^*(x)$ $=$ $\theta(x)x-\Gamma(\theta(x))$.   
The key step is now to introduce the new probability distribution $\mu_\theta$ as in \reff{mutheta} and $\P_\theta$ the corresponding probability measure on $(\Omega,\Fc)$ with likelihood ratio~:
\beqs
\frac{d\P_\theta}{d\P} &=& \Prod_{i=1}^n \frac{d\mu_\theta}{d\mu}(X_i) \; = \; \exp\Big(\theta S_n - n \Gamma(\theta)\Big). 
\enqs
Then, we have by \reff{changeesp} and for all $\eps$ $>$ $0$~: 
\beqs
\P \Big[\frac{S_n}{n} \in  [x,x+\eps) \Big]   &=& \E_\theta \Big[ \exp\Big(- \theta S_n + n \Gamma(\theta)\Big) 1_{\frac{S_n}{n}\in [x,x+\eps)} \Big] \\
&=&  e^{-n(\theta x - \Gamma(\theta))}   \E_\theta \Big[ \exp\Big(- n \theta  ( \frac{S_n}{n} - x) \Big) 1_{\frac{S_n}{n}\in [x,x+\eps)} \Big]  \\
& \geq &  e^{-n(\theta x - \Gamma(\theta))}  e^{-n|\theta|\eps} \P_\theta\Big[ \frac{S_n}{n} \in  [x,x+\eps) \Big],
\enqs
and so
\beq \label{interlower}
\frac{1}{n} \ln \P \Big[\frac{S_n}{n} \in  [x,x+\eps) \Big]  &\geq & - [\theta x - \Gamma(\theta)]  -  |\theta|\eps  
+ \frac{1}{n} \ln \P_\theta\Big[ \frac{S_n}{n} \in  [x,x+\eps) \Big].  
\enq
Now, since $\Gamma'(\theta)$ $=$ $x$, we have $\E_\theta[X_1]$ $=$ $x$, and by the law of large numbers and CLT~: 
$\lim_n  \P_\theta\Big[ \frac{S_n}{n} \in  [x,x+\eps) \Big]$ $=$ $1/2$ $(>0)$.  We also have $\Gamma^*(x)$ $=$ $\theta x -\Gamma(\theta)$.  
Therefore, by sending $n$ to infinity and then $\eps$ to zero in \reff{interlower}, we get \reff{interlowery}. 

Now, if the supremum in $\Gamma^*(x)$ is not attained, we can find a sequence $(\theta_k)_k$ $\nearrow$ $\infty$, such  that 
$\theta_k x  - \Gamma(\theta_k)$ $\rightarrow$ $\Gamma^*(x)$. Since  $\E[e^{\theta_k(X_1-x)}1_{X_1<x}]$ $\rightarrow$ $0$, we then get 
\beqs
\E[e^{\theta_k(X_1-x)}1_{X_1\geq x}] & \rightarrow & e^{-\Gamma^*(x)}, 
\enqs
as $k$ goes to infinity.  This is possible only if $\P[X_1 > x]$ $=$ $0$ and $\P[X_1=x]$ $=$ $e^{-\Gamma^*(x)}$.  By the i.i.d. assumption on the $X_i$, 
this implies $\P[S_n/n \geq x]$ $\geq$ $(\P[X_1\geq x])^n$ $=$ $e^{-n\Gamma^*(x)}$, which proves \reff{interlowery}.

Suppose now that $\mu$ is of unbounded support, and  fix  $M$ large enough 
s.t. $\mu([-M,M])$ $>$ $0$. By the preceding proof,  the lower bound \reff{interlowery} holds with the law of $S_n/n$ conditional 
on $\{|X_i|\leq M, i=1,\ldots,n\}$, and with a c.g.f. equal to the c.g.f. of  the conditional law of $X_1$ given $|X_1|$ $\leq$ $M$~: 
\beq 
& & \lim_{\eps\rightarrow 0} \liminf_{n\rightarrow\infty} \frac{1}{n} 
 \ln \P \Big[\frac{S_n}{n} \in  [x,x+\eps) \Big|  |X_i| \leq M, i=1,\ldots,n \Big]  \nonumber \\
& \geq & - \tilde \Gamma_M^*(x) \; := \; - \sup_{\theta\in\R}[ \theta x - \tilde \Gamma_M(\theta)],  \label{interloweryM}
\enq
with  $\tilde \Gamma_M(\theta)$ $=$ $\ln \E[e^{\theta X_1}| |X_1| \leq M]$ $=$ $\Gamma_M(\theta)-\ln \mu([-M,M])$, $\Gamma_M(\theta)$ 
$=$ $\ln \E[e^{\theta X_1} 1_{|X_1|\leq M}]$.  Now, by writing from Bayes formula 
that $\P\Big[\frac{S_n}{n} \in  [x,x+\eps) \Big]$ $=$ $\P \Big[\frac{S_n}{n} \in  [x,x+\eps) \Big|  |X_i| \leq M, i=1,\ldots,n \Big]$.
$(\mu([-M,M]))^n$, we get with \reff{interloweryM}
\beqs
& & \lim_{\eps\rightarrow 0} \liminf_{n\rightarrow\infty} \frac{1}{n} \ln \P \Big[\frac{S_n}{n} \in  [x,x+\eps) \Big]  \\
& \geq &  \lim_{\eps\rightarrow 0} \liminf_{n\rightarrow\infty} \frac{1}{n} 
 \ln \P \Big[\frac{S_n}{n} \in  [x,x+\eps) \Big|  |X_i| \leq M, i=1,\ldots,n \Big]  + \ln \mu([-M,M]) \\
 & \geq &  -  \Gamma_M^*(x) \; := \; - \sup_{\theta\in\R}[ \theta x -  \Gamma_M(\theta)]. 
\enqs
The required result is obtained by sending $M$ to infinity.  Notice also finally that $\inf_{y\geq x} \Gamma^*(y)$ $=$ $\Gamma^*(x)$ since $\Gamma^*$ is nondecreasing on $[\E X_1,\infty)$. 
\ep

\vspace{2mm}

\noindent {\bf Examples}

\noindent 1) Bernoulli distribution~: for $X_1$ $\sim$ $\Bc(p)$, we have $\Gamma^*(x)$ $=$   $x\ln\big(\frac{x}{p}\big) 
+ (1-x)\ln\big(\frac{1-x}{1-p}\big)$ for $x$ $\in$ $[0,1]$ and $\infty$ otherwise.

\noindent 2) Poisson distribution~: for $X_1$ $\sim$ $\Pc(\lambda)$, we have $\Gamma^*(x)$ $=$   $x\ln\big(\frac{x}{\lambda}\big) + \lambda - x$ for 
$x$ $\geq$ $0$ and $\infty$ otherwise.

\noindent 3) Normal distribution~:  for $X_1$ $\sim$ $\Nc(0,\sigma^2)$, we have $\Gamma^*(x)$ $=$ $\frac{x^2}{2\sigma^2}$, $x$ $\in$ $\R$.  

\noindent 2) Exponential distribution~: for $X_1$ $\sim$ $\Ec(\lambda)$, we have $\Gamma^*(x)$ $=$ $\lambda x - 1 - \ln(\lambda x)$ for $x$ $>$ $0$ and $\Gamma^*(x)$ $=$ $\infty$ otherwise.

\begin{Remark}   
 {\rm   Cramer's theorem possesses a multivariate counterpart dealing with the large deviations of the empirical means of i.i.d. random vectors in $\R^d$.  
 }
\end{Remark}

\begin{Remark} \label{remIS}
{\rm (Relation with importance sampling)

\noindent  Fix $n$ and let us consider the estimation of $p_n$ $=$ $\P[S_n/n \geq x]$. A standard estimator for $p_n$ is the average with $N$ independent copies of $X$ $=$ $1_{S_n/n\geq x}$. However, as shown in the introduction,  for large $n$, $p_n$ is small, and the relative error of this 
estimator is  large.  By using an exponential change of measure $\P_\theta$ with likelihood ratio
\beqs
\frac{d\P_\theta}{d\P} &=& \exp\big( \theta S_n - n \Gamma(\theta)\big), 
\enqs
so that 
\beqs
p_n &=& \E_\theta \Big[ \exp\big( - \theta S_n + n \Gamma(\theta)\big) 1_{\frac{S_n}{n}\geq x} \Big],
\enqs 
we have an importance sampling (IS)  (unbiased) estimator of $p_n$, by taking  the average of independent replications of 
\beqs
 \exp\big( - \theta S_n + n \Gamma(\theta)\big) 1_{\frac{S_n}{n}\geq x}. 
\enqs
The parameter  $\theta$ is  chosen in order to minimize the variance of this estimator, or equivalently its second moment~: 
\beq
M_n^2(\theta,x) & = & \E_\theta\Big[  \exp\big( - 2\theta S_n + 2n \Gamma(\theta)\big) 1_{\frac{S_n}{n}\geq x} \Big]  \nonumber \\
&\leq &   \exp\big( - 2n(\theta  x  - \Gamma(\theta))\big)  \label{optim2}
\enq
By noting from Cauchy-Schwarz's inequality that $M_n^2(\theta,x)$ $\geq$ $p_n^2$ $=$ $\P[S_n/n \geq x ]$ $\simeq$ $C e^{-2n\Gamma^*(x)}$ as $n$ goes to infinity,  from Cramer's theorem, we see that the fastest possible exponential rate of decay of $M_n^2(\theta,x)$ is twice the rate of 
the probability itself, i.e. $2\Gamma^*(x)$. Hence, from \reff{optim2}, and with the choice of $\theta$ $=$ $\theta_x$ s.t. $\Gamma^*(x)$ $=$ $\theta_x x - \Gamma(\theta_x)$, we get an 
asymptotic optimal IS estimator in the sense that~:
\beqs
\lim_{n\rightarrow\infty} \frac{1}{n} \ln M_n^2(\theta_x,x) &=& 2 \lim_{n\rightarrow\infty} \frac{1}{n} \ln p_n. 
\enqs
This  parameter $\theta_x$  is such that $\E_{\theta_x}[S_n/n]$ $=$ $x$ so that the event $\{S_n/n \geq x\}$  is no more rare 
under $\P_{\theta_x}$, and is precisely the parameter used in the derivation of the large deviations result in Cramer's theorem.  
}
\end{Remark}

\subsection{Some general  principles in large deviations}  \label{paraprin}

In this section, we give  some general  principles in large deviations theory.  We refer to the classical references \cite{demzei98} or \cite{dupell97}  for a 
detailed treatment  on the subject.

We first give the formal definition of a large deviation principle (LDP).  Consider a sequence $\{Z^\eps\}_\eps$ on $(\Omega,\Fc,\P)$ valued in some topological space $\Xc$.  The LDP characterizes the limiting behaviour as $\eps$ $\rightarrow$ $0$ of the family of probability measures 
$\{\P[ Z^\eps \in dx]\}_\eps$ on 
$\Xc$ in terms of a {\it rate function}.  A  rate function $I$ is a lower semicontinuous function mapping $I$ $:$ $\Xc$ $\rightarrow$ $[0,\infty]$. It is a 
{\it good} rate function if the level sets $\{ x \in \Xc~: I(x) \leq M\}$ are compact for all $M$ $<$ $\infty$.  

The sequence $\{Z^\eps\}_\eps$ satisfies a LDP on $\Xc$ with rate function $I$ (and speed $\eps$) if~:

\noindent (i) {\it Upper bound}~: for any closed subset $F$ of $\Xc$
\beqs
\limsup_{\eps\rightarrow 0} \eps \ln \P[Z^\eps \in F ] & \leq & - \inf_{x \in F} I(x). 
\enqs
\noindent (ii) {\it Lower bound}~: for any open subset $G$ of $\Xc$
\beqs
\liminf_{\eps\rightarrow 0} \eps \ln \P[Z^\eps \in G ] & \geq & - \inf_{x \in G} I(x). 
\enqs
If $F$ is a subset of $\Xc$ s.t. $\inf_{x\in F^o}I(x)$ $=$ $\inf_{x\in \bar F}I(x)$ $:=$ $I_F$, then 
\beqs
\lim_{\eps\rightarrow 0}  \eps \ln \P[Z^\eps \in F ]  &=& - I_F,
\enqs
which formally means that $\P[Z^\eps \in F]$ $\simeq$ $C e^{-I_F/\eps}$ for some constant $C$.  The classical Cramer's theorem considered the case of 
the empirical mean  $Z^\eps$ $=$ $S_n/n$ of  i.i.d. random variables in $\R^d$,  with $\eps$ $=$ $1/n$.   Further main results in large deviations theory are the G\"artner-Ellis theorem, which is a version of Cramer's theorem  where independence is weakened to the existence of 
\beqs
\Gamma(\theta) &:=& \lim_{\eps\rightarrow 0} \eps \ln \E\big[e^{\frac{\theta .Z^\eps}{\eps}}\big], \;\; \theta \in \R^d. 
\enqs
LDP is then stated for the sequence $\{Z^\eps\}_\eps$ with a rate function equal to the Fenchel-Legendre transform of $\Gamma$~: 
\beqs
\Gamma^*(x) &=& \sup_{\theta\in\R^d} [\theta.x - \Gamma(\theta)], \;\;\; x \in \R^d. 
\enqs
Other results in large deviations theory include Sanov's theorem, which gives rare events asymptotics for empirical distributions. In many problems, the interest is in rare events that depend on random process, and the corresponding asymptotics probabilities, usually called sample path large deviations, 
were developed by  Freidlin-Wentzell and Donsker-Varadhan.  
For instance, the problem of diffusion exit from a domain is an important application of Freidlin-Wentzell theory, and occurs naturally in finance, see 
Section \ref{secISlar}.  We briefly summarize these results.  Let $\eps$ $>$ $0$ a (small) positive parameter and consider the stochastic differential equation in $\R^d$ on some interval $[0,T]$,  
\beq \label{eqXeps}
dX_s^\eps &=& b_\eps(s,X_s^\eps) ds + \sqrt{\eps} \sigma(s,X_s^\eps)dW_s, 
\enq
and suppose that  there exists a Lipschitz function $b$ on $[0,T]\times\R^d$ s.t. 
\beqs
\lim_{\eps\rightarrow 0}  b_\eps &=& b, 
\enqs 
uniformly  on compact sets.  Given an open set $\Gamma$ of $[0,T]\times\R^d$, we consider the exit time from $\Gamma$, 
\beqs
\tau^\eps_{_{t,x}} &=& \inf\big\{ s \geq t~: X_s^{\eps,t,x} \notin \Gamma \big\},
\enqs
and the  corresponding exit probability
\beqs
v_\eps(t,x) &=&  \P[ \tau^\eps_{_{t,x}} \leq T], \;\;\; (t,x) \in [0,T]\times\R^d. 
\enqs
Here $X^{\eps,t,x}$ denotes the solution to \reff{eqXeps} starting from $x$ at time $t$.  It is well-known that the process $X^{\eps,t,x}$ converge 
to $X^{0,t,x}$ the solution to the ordinary differential equation
\beqs
dX_s^0 &=& b(s,X_s^0) ds, \;\;\;\;\;  X_t^0 = x. 
\enqs 
In order to ensure that  $v_\eps$ goes to zero, we assume that for all $t$ $\in$ $[0,T]$, 
\beqs
{\bf (H)} \;\;\;\;\;\;\;  x \in \Gamma & \Longrightarrow & X_s^{0,t,x}  \in \Gamma, \;\;\; \forall s \in [t,T]. 
\enqs
Indeed, under {\bf (H)},   the system \reff{eqXeps} tends, when $\eps$ is small,  to stay inside $\Gamma$, so that the event 
$\{\tau^\eps_{_{t,x}} \leq T\}$  is rare.  The large deviations asymptotics of $v_\eps(t,x)$, when $\eps$ goes to zero,  was initiated  by Varadhan and Freidlin-Wentzell by probabilistic arguments.  An alternative approach, introduced by Fleming, connects  this theory with optimal control and Bellman equation, and is developed within the 
theory of viscosity solutions,  see e.g. \cite{bar94}.  We sketch here this approach.  
It is well-known that the function $v_\eps$ satisfies the linear PDE
\beq \label{edpueps}
\Dt{v_\eps} +   b_\eps(t,x). D_x v_\eps +  \frac{\eps}{2}  {\rm tr}(\sigma\sigma'(t,x)D_x^2 v_\eps) &=& 0, \;\;\; (t,x) \in [0,T)\times\Gamma
\enq
together with the boundary conditions
\beq
v_\eps(t,x) &=& 1, \;\;\; (t,x) \in [0,T)\times  \partial\Gamma \label{boundueps} \\
v_\eps(T,x) &=& 0, \;\;\; x \in \Gamma.  \label{termueps}
\enq
Here $\partial \Gamma$ is the boundary of $\Gamma$.  
We now make the logarithm transformation 
\beqs
V_\eps &=& - \eps \ln v_\eps.
\enqs
Then, after some straightforward derivation,  \reff{edpueps}  becomes the nonlinear PDE
\beq
- \Dt{V_\eps}  -  b_\eps(t,x). D_x V_\eps  -  \frac{\eps}{2}  {\rm tr}(\sigma\sigma'(t,x)D_x^2 V_\eps) & &  \nonumber \\
 + \; \frac{1}{2} (D_x V_\eps)'  \sigma\sigma'(t,x) D_x V_\eps & = &  0, \;\;\;  (t,x) \in [0,T)\times\Gamma, \label{edpUeps} 
\enq
and the boundary data \reff{boundueps}-\reff{termueps} become
\beq
V_\eps(t,x) &=& 0, \;\;\; (t,x) \in [0,T)\times  \partial\Gamma \label{boundUeps} \\
V_\eps(T,x) &=& \infty, \;\;\; x \in \Gamma.  \label{termUeps}
\enq 
At the limit $\eps$ $=$ $0$,  the PDE \reff{edpUeps} becomes a first-order PDE
\beq
- \Dt{V_0}  -  b(t,x). D_x V_0  +  \frac{1}{2} (D_x V_0)'  \sigma\sigma'(t,x) D_x V_0 & = &  0, \;\;\;  (t,x) \in [0,T)\times\Gamma, \label{edpU0} 
\enq 
with the boundary data \reff{boundUeps}-\reff{termUeps}.   By PDE-viscosity solutions methods and comparison results,  we can prove (see e.g. 
\cite{bar94} or \cite{fleson94}) that $V_\eps$ converges  uniformly on compact subsets of $[0,T)\times\Gamma$, as $\eps$ goes to zero,  to $V_0$ the unique viscosity solution to \reff{edpU0} with the boundary data \reff{boundUeps}-\reff{termUeps}.  Moreover, $V_0$ has a representation in terms of 
control problem.  Consider the Hamiltonian function
\beqs
H(t,x,p) &=&  -b(t,x).p +  \frac{1}{2} p' \sigma\sigma'(t,x) p, \;\;\; (t,x,p) \in [0,T]\times\Gamma \times \R^d,
\enqs
 which is quadratic and in particular convex in $p$. Then, using the Legendre transform, we may rewrite
 \beqs
 H(t,x,p) &=&  \sup_{q \in \R^d} \big[  -  q.p  -   H^*(t,x,q)  \big], 
 \enqs
 where 
 \beqs
 H^*(t,x,q) &=&  \sup_{p \in \R^d} \big[  -  p.q -   H(t,x,p)  \big]  \\
 &=&    \frac{1}{2} (q-b(t,x))' (\sigma\sigma'(t,x))^{-1} (q-b(t,x)), \;\; (t,x,q) \in [0,T]\times\Gamma \times \R^d. 
 \enqs
Hence, the PDE \reff{edpU0} is rewritten as 
\beqs
 \Dt{V_0}  +   \inf_{q \in \R^d} \big[    q. D_x V_0  +   H^*(t,x,q)  \big] &=& 0,  \;\;\;  (t,x) \in [0,T)\times\Gamma,
\enqs
which, together with the boundary data \reff{boundUeps}-\reff{termUeps}, is associated to the value function for the  following   calculus of variations problem~:  for an absolutely continuous  function $x(.)$ on $[0,T)$ and valued in $\R^d$, i.e. $x$ $\in$ $H^1_{loc}([0,T],\R^d)$, we denote $\dot{x}(u)$ $=$ $q_u$  its time derivative, and 
$\tau(x)$ the exit time of $x(.)$ from $\Gamma$. Then, 
\beqs
V_0(t,x) &=& \inf_{x(.) \in \Ac(t,x)}  \int_t^T   H^*(u,x(u),\dot{x}(u)) du,  \;\;\;  (t,x) \in [0,T)\times\Gamma,  \\
&=&  \inf_{x(.) \in \Ac(t,x)}  \int_t^T  \frac{1}{2} (\dot x(u)-b(u,x(u)))' (\sigma\sigma'(u,x(u)))^{-1} (\dot x(u)-b(u,x(u))) du
\enqs
where 
\beqs
\Ac(t,x) &=& \big\{ x(.) \in H^1_{loc}([0,T],\R^d)~:   \;  x(t) = x \; \mbox{ and } \; \tau(x) \leq  T \big\}. 
\enqs
The large deviations result is then stated as
\beq \label{larfrei}
\lim_{\eps\rightarrow 0} \eps \ln v_\eps(t,x) &=& - V_0(t,x),
\enq
and the above limit holds uniformly on compact subsets of $[0,T)\times\Gamma$.   A more precise result may be obtained, which allows to remove 
the above log estimate. This  type of result is developed  in  \cite{flejam92}, and is called sharp large deviations estimate.  It states asymptotic expansion 
(in $\eps)$ of the exit probability for points $(t,x)$ belonging to a set $N$ of 
$[0,T']\times\Gamma$ for some $T'$ $<$ $T$, open in the relative topology,  and s.t.  $V_0$ $\in$ $C^\infty(N)$.   Then, under the condition that 
\beqs
b_\eps &=& b + \eps b_1 +  0(\eps^2), 
\enqs
one has 
\beqs
v_\eps(t,x) &=& \exp\big( - \frac{V_0(t,x)}{\eps} - w(t,x)  \big)(1 + O(\eps)),
\enqs
uniformly on compact sets of $N$, where $w$ is solution to the PDE problem
\beqs
- \Dt{W} - (b-\sigma\sigma' D_xV_0). D_x w &=& \frac{1}{2} {\rm tr}(\sigma\sigma'D_x^2V_0) + b_1 . D_x V_0 \;\;\;\; \mbox{ in } \;\; N \\
w(t,x) &=& 0 \;\;\; \mbox{ on } \;\;\; \Big([0,T)\times\partial\Gamma\Big)  \cup \bar N. 
\enqs
The function $w$ may be represented as 
\beqs
w(t,x) &=& \int_t^\rho \Big(\frac{1}{2} {\rm tr}(\sigma\sigma'D_x^2V_0) + b_1 . D_x V_0\Big)(s,\xi(s)) ds,
\enqs
where $\xi$ is the solution to 
\beqs
\dot \xi(s) &=&  (b-\sigma\sigma' D_xV_0)(s,\xi(s)), \;\;\;  \xi(t) \; = \; x, 
\enqs
and $\rho$ is the exit time (after $t$) of $(s,\xi(s))$ from $N$.

We shall develop more in detail  in the next sections some applications of the G\"artner-Ellis and Freidlin-Wentzell theories in finance.

We end this paragraph by stating the important Varadhan's integral formula, which involves the asymptotics behavior of certain expectations. It extends the well-known method of Laplace for studying the asymptotics of certain integrals  on $\R$~: given a continuous function $\varphi$ from $[0,1]$ into $\R$, 
Laplace's method states that 
\beqs
\lim_{n\rightarrow\infty} \frac{1}{n} \ln \int_0^1 e^{n\varphi(x)} dx &=&  \max_{x\in [0,1]} \varphi(x). 
\enqs
Varadhan result's is  formulated  as follows~:

\begin{Theorem} \label{theovaradh} (Varadhan)
\noindent Suppose that $\{Z^\eps\}_\eps$ satisfies a LDP on $\Xc$ with good rate function $I$, and let $\varphi$ $:$ $\Xc$ $\rightarrow$ $\R$ be any 
continuous function s.t. the following moment condition holds for some $\gamma$ $>$ $1$~:
\beqs
\limsup_{\eps\rightarrow 0} \eps \ln \E\big[e^{\gamma \varphi(Z^\eps)/\eps}\big] & < & \infty. 
\enqs
Then, 
\beq \label{varad}
\lim_{\eps\rightarrow 0} \eps \ln \E\big[e^{\varphi(Z^\eps)/\eps}\big] &=& \sup_{x\in\Xc} \big[\varphi(x) - I(x)\big]. 
\enq
\end{Theorem}
{\bf Proof.}
(a) For simplicity, we show the inequality $\leq$ in \reff{varad}  
when $\varphi$ is bounded on $\Xc$.  Hence, there exists $M$ $\in$ $(0,\infty)$ s.t.  $-M$ $\leq$ $\varphi(x)$ $\leq$ 
$M$ for all $x$ $\in$ $\Xc$. For $N$ positive integer,  and $j$ $\in$ $\{1,\ldots,N\}$, we consider the closed subsets of $\Xc$
\beqs
F_{N,j} &=& \big\{ x \in \Xc~:  - M + \frac{2(j-1)M}{N} \leq \varphi(x) \leq -M +  \frac{2jM}{N} \big\}, 
\enqs
so that $\cup_{j=1}^N F_{N,j}$ $=$ $\Xc$.  We then have from the large deviations upper bound on $(Z^\eps)$, 
\beqs
\limsup_{\eps\rightarrow 0} \eps \ln \E\big[e^{\varphi(Z^\eps)/\eps}\big]  & =  & 
\limsup_{\eps\rightarrow 0} \eps \ln \int_{\Xc}  e^{\varphi(Z^\eps)/\eps} \P[Z^\eps \in dx] \\
& \leq &  \limsup_{\eps\rightarrow 0} \eps \ln \big( \sum_{j=1}^N \int_{F_{N,j}}  e^{\varphi(Z^\eps)/\eps} \P[Z^\eps \in dx]  \big) \\
& \leq &   \limsup_{\eps\rightarrow 0} \eps \ln \big( \sum_{j=1}^N e^{(-M+2jM/N)/\eps}\P[Z^\eps \in F_{N,j}] \big) \\
& \leq &   \limsup_{\eps\rightarrow 0} \eps \ln \big( \max_{j=1,\ldots,N} e^{(-M+2jM/N)/\eps}\P[Z^\eps \in F_{N,j}] \big) \\
& \leq &  \max_{j=1,\ldots,N} \Big( -M +  \frac{2jM}{N}  + \limsup_{\eps\rightarrow 0} \eps  \ln \P[Z^\eps \in F_{N,j}]  \Big) \\
& \leq &   \max_{j=1,\ldots,N} \Big( -M +  \frac{2jM}{N} +   \sup_{x\in F_{N,j}}[- I(x)] \Big)\\
& \leq &   \max_{j=1,\ldots,N}  \Big( -M +  \frac{2jM}{N} + \sup_{x\in F_{N,j}}[\varphi(x) - I(x)]  - \inf_{x \in F_{N,j}} \varphi(x)\Big) \\
& \leq &  \sup_{x\in F_{N,j}}[\varphi(x) - I(x)]   +  \frac{2M}{N}.
\enqs
By sending $N$ to infinity, we get the inequality $\leq$ in \reff{varad}. 

\noindent (b) To prove the reverse inequality, we fix an arbitrary point $x_0$ $\in$ $\Xc$, an arbitrary $\delta$ $>$ $0$, and we consider the open set 
$G$ $=$ $\{ x \in \Xc~:$ $\varphi(x) > \varphi(x_0) - \delta\}$.  Then, we have from the large deviations lower bound on $(Z^\eps)$, 
\beqs
\liminf_{\eps\rightarrow 0} \eps \ln \E\big[e^{\varphi(Z^\eps)/\eps}\big]  & \geq & 
\liminf_{\eps\rightarrow 0} \eps \ln \E\big[e^{\varphi(Z^\eps)/\eps}1_{Z^\eps \in G} \big] \\
& \geq &  \varphi(x_0) - \delta +  \liminf_{\eps\rightarrow 0} \eps \ln  \P[Z^\eps \in G]  \\
& \geq & \varphi(x_0) - \delta - \inf_{x\in G} I(x)  \\
& \geq & \varphi(x_0) - I(x_0) - \delta. 
\enqs
Since $x_0$  $\in$ $\Xc$ and $\delta$ $>$ $0$ are arbitrary, we get the required result. 
\ep

\begin{Remark}
{\rm The relation \reff{varad} has the following  interpretation. By writing formally the  LDP for $(Z^\eps)$ with rate function $I$ as 
$\P[Z^\eps \in dx]$ $\simeq$ $e^{-I(x)/\eps}dx$,  we can write 
\beqs
\E\big[e^{\varphi(Z^\eps)/\eps}\big] \; = \; \int e^{\varphi(x)/\eps} \P[Z^\eps \in dx] & \simeq & \int e^{(\varphi(x)-I(x))/\eps}  dx \\
&\simeq & C \exp\Big(\frac{{\sup_{x\in\Xc}(\varphi(x)-I(x))}}{\eps}\Big).  
\enqs
As in Laplace's method, Varadhan's formula states that to exponential order, the main contribution to the integral is due to the largest value of the exponent. 
}
\end{Remark}

\section{Ruin probabilities in  risk theory}

\setcounter{equation}{0} \setcounter{Assumption}{0}
\setcounter{Theorem}{0} \setcounter{Proposition}{0}
\setcounter{Corollary}{0} \setcounter{Lemma}{0}
\setcounter{Definition}{0} \setcounter{Remark}{0}

\subsection{The classical  ruin  problem}

\subsubsection{The insurance model}

We consider an insurance company earning premiums at a constant rate $p$ per unit of time, and paying claims that arrive at the jumps of a Poisson process with intensity $\lambda$.  We denote by $N_t$ the number of claims arriving in $[0,t]$,  by  $T_n$, $n$ $\geq$ $1$, the arrival times of the claim, 
and by $\xi_1$ $=$ $T_1$, $\xi_n$ $=$ $T_n-T_{n-1}$, $n$ $\geq$ $2$,   the interarrival times,  
which are then i.i.d. exponentially distributed with finite mean 
$E \xi_1$ $=$ $1/\lambda$.   The size of the $n$-th claim  is denoted $Y_n$, and we assume that the claim sizes $Y_n$, $n$ $\in$ $\N^*$, are (positive) i.i.d., and independent of the Poisson process.   Starting from an initial reserve $x$ $>$ $0$, the risk reserve process $X_t$ $=$ $X_t^x$,  $t$ $\geq$ $0$, of the insurance company is then given by~: 
\beq \label{defriskX}
X_t^x &=& x + p t - \sum_{i=1}^{N_t} Y_i. 
\enq
The probability of ruin with infinite horizon is 
\beqs
\psi(x) &=& \P[ \tau_x < \infty], 
\enqs
where $\tau_x$ $=$ $\inf\{ t \geq 0~: X_t^x < 0\}$ is the time of ruin.   We are   interested in the estimation of the ruin probability, in particular for large values of the initial reserve.

\subsubsection{The Cramer-Lundberg estimate}

The Cramer-Lundberg  approximation concerns the estimation of the  ruin probability $\psi(x)$, and is one of the most celebrated result of risk theory. 
There are several approaches for deriving such a result. We follow in this paragraph a method based on large deviations arguments and change of probability measures.

First, we easily see,  by the strong law of large numbers, that 
\beqs
\frac{1}{t} \sum_{i=1}^{N_t} Y_i & \rightarrow &  \rho  \;\;\; a.s., \;\;\;\;\; t \rightarrow \infty, 
\enqs 
where $\rho$ $=$ $\lambda \E[Y_1]$ $>$ $0$  is interpreted as the average amount of claim per unit of time.  
The safety loading $\eta$  plays a key role in ruin probability. It is  defined as the relative amount by which the premium rate exceeds $\rho$~: 
\beqs
\eta &=& \frac{p - \rho}{\rho} \;\;\; \Longleftrightarrow \;\;\; p \; = \; (1+\eta) \rho. 
\enqs
Hence $X_t^x/t$ $\rightarrow$ $p-\rho$ $=$ $\rho\eta$ when $t$ goes to infinity. Therefore, if $\eta$ $<$ $0$, $X_t^x$ $\rightarrow$ $-\infty$, and we clearly have $\psi(x)$ $=$ $1$ for all $x$.  For $\eta$ $=$ $0$, we can also show that $\lim\sup  X_t^x$ $=$ $-\infty$ so that $\psi(x)$ $=$ $1$.   In the sequel, we make the net profit assumption~: 
\beq \label{netprofit}
\eta \; = \; \frac{p - \lambda \E[Y_1]}{\lambda \E[Y_1]} \; > \;  0,
\enq
which ensures that  the probability of ruin is less than $1$.  

Since ruin may occur only  at the arrival of a claim, i.e. when $X$ jumps downwards, it suffices to consider the discrete-time process embedded at the jumps of the Poisson process. We then define the discrete-time process  $X_{T_n}^x$, $n$ $\geq$ $1$, so that 
\beqs
\psi(x) &=&  \P[ \sigma_x < \infty],
\enqs
where $\sigma_x$ $=$ $\inf\{ n \geq 1~: X_{T_n}^x < 0\}$ $=$ $\inf\{ n \geq 1~: S_n > x\}$, and $S_n$ $=$ $x- X_{T_n}^x$ is the net payout up to the $n$-th claim and given by the random walk~: 
\beqs
S_n &=&  Z_1 + \ldots + Z_n, \;\;\;\; Z_i \; = \; Y_i - p \xi_i, \;\; i \in \N^*. 
\enqs
The r.v. $Z_i$ are i.i.d. and satisfy under the net profit condition, $\E[Z_1]$ $<$ $0$.  We denote by $\Gamma_Z$ the c.g.f. of the $Z_i$, and we see  that by independence of $Y_i$ and $\xi_i$~:
\beqs
\Gamma_Z(\theta) &=& \Gamma_Y(\theta) + \Gamma_\xi(-p\theta) \\
&=&   \Gamma_Y(\theta) + \ln \Big( \frac{\lambda}{\lambda + p\theta} \Big), \;\;\; \theta > - \frac{\lambda}{p},
\enqs
where $\Gamma_Y$ (resp. $\Gamma_\xi$)  is the c.g.f. of the $Y_i$ (resp. $\xi_i$).   For any $\theta$ in the domain of $\Gamma_Z$, we consider the exponential change of measure with parameter $\theta$, and since $\sigma_x$ is a stopping time in the filtration of $(Z_1,\ldots,Z_n)$, 
we apply formula \reff{changeesptau} to write the ruin probability as an $\E_\theta$ expectation~: 
\beq
\psi(x) &=&   \P[ \sigma_x < \infty] \; = \; \E\big[ 1_{\sigma_x < \infty} \big] \nonumber \\
&=& \E_\theta \Big[ 1_{\sigma_x < \infty}   \exp\big(-\theta S_{\sigma_x} + \sigma_x \Gamma_Z(\theta) \big) \Big]. \label{changeruin}
\enq
We now assume $Y$ has a light-tailed distribution, i.e.~: there exists $\bar\theta$ $\in$ $(0,\infty]$ s.t. $\Gamma_y(\theta)$ $<$ $\infty$ for $\theta$ $<$ 
$\bar\theta$, and $\Gamma_Y(\theta)$ $\rightarrow$ $\infty$ as $\theta$ $\nearrow$ $\bar\theta$.   In this case,  the c.g.f.  $\Gamma_Z$ of the $Z_i$ is  
finite on $(-\lambda/p,\bar\theta)$,  it  is differentiable in $0$ with $\Gamma_Z'(0)$ $=$ $\E[Z_1]$ $<$ $0$  under 
the net profit condition \reff{netprofit}.  Moreover, since $\E[Y_1]$ $>$ $0$ and $Y_1$ is independent of $\xi_1$, we see that $\P[Z_1 > 0]$ $>$ $0$, which implies that $\Gamma_Z(\theta)$ goes to infinity as $\theta$ goes to $\bar\theta$.  
By convexity of $\Gamma_Z$ and recalling that $\Gamma_Z(0)$ $=$ $0$, 
we deduce the existence of an unique $\theta_{_{L}}$ $>$ $0$ s.t. $\Gamma_Z(\theta_{_{L}})$ $=$ $0$. This unique positive $\theta_{_{L}}$ is the solution to the so-called  {\it Cramer-Lundberg} equation~: 
\beqs
\Gamma_Y(\theta_{_{L}}) +  \ln \Big( \frac{\lambda}{\lambda + p\theta_{_{L}}} \Big) &=& 0,
\enqs
which is also written equivalently in~: 
\beq \label{eqgammalundberg}
\gamma_{_Y}(\theta_{_{L}}) &=& \frac{p\theta_{_{L}}}{\lambda},
\enq
where $\gamma_{_Y}$ $=$ $\exp(\Gamma_Y(\theta))-1$ $=$ $\int e^{\theta y} F_Y(dy)-1$ is the shifted ($\gamma_Y(0)$ $=$ $0$) 
moment generating function of $Y_i$, and 
$F_Y$ is the distribution function of the claim sizes $Y_i$.   $\theta_{_{L}}$ is called {\it adjustment coefficient} (or sometimes {\it Lundberg exponent}). 
Notice also  that  by convexity of $\Gamma_Z$, we have $\Gamma_Z'(\theta_{_{L}})$ 
$>$ $0$.  Hence,  under $\P_{\theta_{_{L}}}$, the random walk  has positive drift $\E_{\theta_{_{L}}}[Z_n]$ $=$ $\Gamma_Z'(\theta_{_{L}})$ $>$ $0$, and this implies 
$\P_{\theta_{_{L}}}[\sigma_x < \infty]$ $=$ $1$.  For this choice of $\theta$ $=$ $\theta_{_{L}}$, \reff{changeruin} becomes
\beq
\psi(x) &=&  \E_{\theta_{_{L}}}\Big[ e^{-\theta_{_{L}} S_{\sigma_x}}  \Big] \; = \; e^{-\theta_{_{L}} x}  \E_{\theta_{_{L}}}\Big[ e^{-\theta_{_{L}}
(S_{\sigma_x} -x)}  \Big].  
\enq
By noting that the overshoot $S_{\sigma_x} -x$ is nonnegative, we obtain the Lundberg's inequality on the ruin probability~: 
\beq
\psi(x) & \leq &  e^{-\theta_{_{L}} x} , \;\;\; \forall x > 0. 
\enq
Moreover, by renewal's theory, the overshoot $R^x$ $=$ $S_{\sigma_x}-x$ has a limit $R^\infty$ (in the sense of weak convergence with respect to 
$\P_{\theta_{_{L}}}$), when $x$ goes to infinity, and therefore $\E_{\theta_{_{L}}}[e^{-\theta_{_{L}} R^x}]$ converges to some positive constant $C$.  We then get the classical approximation for large values of the initial reserve~: 
\beqs
\psi(x) & \simeq & C e^{-\theta_{_{L}} x},
\enqs
as $x$ $\rightarrow$ $\infty$, which  implies  a large deviation type estimation 
\beq \label{estilund}
\lim_{x\rightarrow\infty} \frac{1}{x} \ln \psi(x) &=& - \theta_{_{L}}. 
\enq
Further details and extensions can be found in \cite{asm00} or \cite{embklumik97}.  They concern more general processes (e.g. Levy proceses) for the risk reserve,  heavy-tailed distribution for the claim size ...  In the next paragraph, we study an extension of the classical ruin problem, developed by 
\cite{gaigrasch03} and \cite{hipsch04},  where the insurer has the additional opportunity to invest in a stock market.

\vspace{2mm}

\noindent {\bf Application~:  Importance sampling for the ruin probability estimation}

\vspace{1mm}

\noindent  From  the perspective of  estimation of the ruin probability $\psi(x)$, and by choosing the Lundberg exponent $\theta_L$,  we have an unbiased estimator with the associated importance sampling  estimator based on Monte-Carlo simulations of 
\beqs
\psi(x) &=& \E_{\theta_L}\big[e^{-\theta_L S_{\sigma_x}}\big]. 
\enqs
 Since, obviously, $S_{\sigma_x}$ $>$ $x$, the second order moment  of this estimator  satisfies
 \beqs
 M^2(\theta_L,x) & = &   \E_{\theta_L}\big[e^{-2\theta_L S_{\sigma_x}}\big] \; \leq \;  e^{-2\theta_L x}. 
 \enqs
 On the other hand, by Cauchy-Schwarz's inequality, the second moment of any unbiased estimator must be as large as the square of the ruin probability, and we have seen that this probability is $O(e^{-\theta_L x})$.   Therefore, the IS estimator based on $\theta_L$ is asymptotically optimal as $x$ 
 $\rightarrow$  $\infty$~: 
 \beqs
 \lim_{x\rightarrow\infty}  \frac{1}{x} \ln M^2(\theta_L,x) &=& 2  \lim_{x\rightarrow\infty}  \frac{1}{x} \ln \psi(x) \; ( = \; 2 \theta_L).  
 \enqs

\begin{Remark}
{\rm  For  any $x$ $>$ $0$, $\theta$ $>$ $0$, consider the process
\beqs
M_t(x,\theta) &=& \exp(-\theta X_t^x), \;\;\; t \geq 0. 
\enqs
where $X_t^x$ is the risk reserve process defined in \reff{defriskX}. A standard calculation shows that for all $t$ $\geq$ $0$, 
\beqs
\E\big[ M_t(0,\theta) \big] &=&  e^{(\lambda \gamma_{_Y}(\theta) - p\theta)t}.
\enqs
Moreover, since $X_t^x$  has stationary independent increments, and by denoting $\F$ $=$ $(\Fc_t)_{t\geq 0}$ the filtration generated by  the risk 
reserve $X$,  we have  for all $0\leq t\leq T$
\beq
\E\big[  M_T(x,\theta) \big|  \Fc_t \big] &=&  \E\big[  e^{-\theta X_T^x}  \big|  \Fc_t \big]  \nonumber \\
&=&  M_t(x,\theta)  \E\big[  e^{-\theta (X_T^x - X_t^x) }  \big|  \Fc_t \big]  \nonumber \\
&=&  M_t(x,\theta)  \E\big[  M_{T-t}(0,\theta)   \big]  \nonumber \\
&=&  M_t(x,\theta) \; e^{(\lambda \gamma_{_Y}(\theta) - p\theta)(T-t)}. \label{martin}
\enq
Hence, for the choice of $\theta$ $=$ $\theta_{_{L}}$~: the adjustment coefficient, the process $M_t(x,\theta_{_{L}})$, $t\geq 0$, is a $(\P,\F)$-martingale.  
The use of this martingale property  in the derivation of ruin estimate was initiated by Gerber \cite{ger73}.  We show in the next  paragraph how to extend 
this idea for a ruin problem with investment in a stock market.   
}
\end{Remark}

\subsection{Ruin probabilities and optimal investment}

\subsubsection{The insurance-finance model}

In the setting of the classical model described in the previous paragraph, we consider  the additional feature that the insurance company is also 
allowed to invest in some stock market, modeled by a geometric brownian motion~: 
\beqs
dS_t &=& b S_t dt  + \sigma S_t dW_t, 
\enqs
where $b,\sigma$ are constants, $\sigma$ $>$ $0$, and $W$ is a standard brownian motion, independent of the risk  reserve $X$ as defined in 
\reff{defriskX}.  We denote by $\F$ $=$ $(\Fc_t)_{t\geq 0}$ the filtration generated by $X$ and $S$.  The insurer may invest at any time $t$ an amount of money $\alpha_t$ in the stock, and the rest in the bond (which in the present model yields no interest). The set $\Ac$  of admissible investment strategies is  defined as the set of $\F$-adapted processes $\alpha$ $=$ $(\alpha_t)$ s.t. $\int_0^t \alpha_s^2 ds$ $<$ $\infty$ a.s. Given an initial capital $x$ $\geq$ $0$, and  an admissible investment control $\alpha$, the insurer's  wealth process can then be written as  
\beqs
V_t^{x,\alpha} &=&  X_t^x  + \int_0^t \frac{\alpha_u}{S_u} dS_u  \\
&=&  x + p t - \sum_{i=1}^{N_t} Y_i + \int_0^t \alpha_u ( bdu + \sigma dW_u), \;\;\;  t \geq 0. 
\enqs
We define  the infinite time ruin probability
\beqs
\psi(x,\alpha) &=& \P[ \tau_{x,\alpha} < \infty],
\enqs
where  $\tau_{x,\alpha}$ $=$ $\inf\{ t \geq 0~: V_t^{x,\alpha}  < 0\}$ is the time of ruin, depending on the initial wealth $x$ and the investment strategy 
$\alpha$.   We are interested in  the minimal ruin probability of the insurer
\beqs
\psi^*(x) &=& \inf_{\alpha\in\Ac} \psi(x,\alpha).
\enqs

\subsubsection{Asymptotic ruin probability estimate}

The main result is an asymptotic large deviation estimation for the minimal ruin probability when the initial reserve goes to infinity~:

\begin{Theorem} \label{theoinsurfin}
We have
\beq \label{estifininsur}
\lim_{x\rightarrow\infty} \frac{1}{x} \ln \psi^*(x) &=& - \theta^*,
\enq
where $\theta^*$ $>$ $0$ is he unique positive solution to the equation
\beq \label{eqgamma}
\gamma_{_Y}(\theta) &=& p \frac{\theta}{\lambda} + \frac{b^2}{2\sigma^2\lambda}. 
\enq
Here $\gamma_{_Y}(\theta)$ $=$ $\E[e^{\theta Y_1}]-1$ is the shifted moment generating function of the claim size. 
Moreover,  the constant strategy $\alpha^*$ $=$ $\frac{b}{\sigma^2\theta^*}$ 
is asymptotically optimal in the sense that
\beqs
\lim_{x\rightarrow\infty} \frac{1}{x} \ln \psi(x,\alpha^*) &=& - \theta^*. 
\enqs
Finally, if  $b$ $\neq$ $0$,  $\theta^*$ $>$ $\theta_{_L}$ the Lundberg exponent.  
\end{Theorem}

\begin{Remark}
{\rm  The estimate \reff{estifininsur} is analogue to the classical Lundberg estimate \reff{estilund}  without investment. The exponent is larger than the Lundberg one, and thus one gets a sharper bound on the minimal ruin probability. Moreover, the trading strategy yielding 
the optimal asymptotic exponential decay consists in holding a fixed (explicit) amount in the risky asset.  This surprising result, in apparent contradiction with the common believe that `rich' companies should invest more than `poor' ones, is explained by the fact that minimization of ruin probability is an extremely conservative criterion. 
}
\end{Remark}

We  follow the martingale approach of Gerber for stating  this theorem. We emphasize the main steps of the proof.  Let us introduce, for fixed 
$x,\theta$ $\in$ $\R_+$, and $\alpha$ $\in$ $\Ac$, the process
\beqs
M_t(x,\theta,\alpha) &=& \exp(-\theta V_t^{x,\alpha}), \;\;\; t \geq 0. 
\enqs
Then, a straightforward calculation shows that for any constant process $\alpha$ $=$ $a$, and $t$ $\geq$ $0$,
\beqs
\E\big[M_t(0,\theta,a)\big] &=& \E\big[e^{-\theta(pt - \sum_{i=1}^{N_t} Y_i + abt + a\sigma W_t)}\big] \\
&=&   e^{-\theta(p+ab)t} \E\big[e^{\theta \sum_{i=1}^{N_t} Y_i}\big] \E\big[e^{-\theta a \sigma W_t}\big] \\
&=&    e^{-\theta(p+ab)t} e^{\gamma_{_Y}(\theta)\lambda t} e^{\frac{\theta^2a^2\sigma^2}{2}t} \\
&=&   e^{f(\theta,a)t},
\enqs
where 
\beqs
f(\theta,a) &=& \lambda \gamma_{_Y}(\theta) -  p\theta  - a b\theta + \frac{1}{2} a^2 \theta^2 \sigma^2. 
\enqs
We recall that under the assumption of light-tailed distribution on the claim size $Y_i$, the shifted moment  generating function $\gamma_{_Y}$ is 
finite  and convex on $(-\infty,\bar\theta)$ for  some $\bar\theta$ $\in$ $(0,\infty]$,   and $\gamma_{_Y}$ $\rightarrow$ $\infty$ as $\theta$ goes to 
$\bar\theta$. Moreover, recalling that  $\E[Y_1]$ $\geq$ $0$, then $\gamma_{_Y}$ is increasing on $[0,\bar\theta)$ since $\gamma'_{_Y}(\theta)$ $=$ 
$\E[Y_1e^{\theta Y_1}]$ $>$ $\E[Y_1]$ for $0$ $<$ $\theta$ $<$ $\bar\theta$. 
Now,  we see that for all $\theta$ $>$ $0$, 
\beqs
\bar f(\theta) \; := \; \inf_{a\in\R} f(\theta,a) &=&  \lambda \gamma_{_Y}(\theta) -  p\theta  -  \frac{b^2}{2\sigma^2},
\enqs
with an infimum attained for $\hat a(\theta)$ $=$ $b/(\theta\sigma^2)$. From the properties of $\gamma_{_Y}$, we clearly have the existence and uniqueness of $\theta^*$ solution to $\bar f(\theta^*)$ $=$ $0$, i.e.  \reff{eqgamma}.  Since the r.h.s. of \reff{eqgamma} is  just the r.h.s. of \reff{eqgammalundberg}, but shifted by the positive constant $b^2/2\sigma^2\lambda$ (if $b$ $\neq$ $0$), it is also obvious that $\theta^*$ $>$ $\theta_{_L}$. 
By choosing $\alpha^*$ $=$ $\hat a(\theta^*)$ $=$ $b^2/(\theta^*\sigma^2)$,  we then have 
\beqs
\bar f(\theta^*) \; = \;  f(\theta^*,\alpha^*) &=& 0.
\enqs
A straightforward calculation also shows that for all $a$ $\in$ $\R$, 
\beq \label{fpos}
f(\theta^*,a) &=& \frac{1}{2} (\theta^*)^2\sigma^2(a - \alpha^*)^2 \; \geq \; 0. 
\enq
Hence, since $V_t^{x,\alpha^*}$ has independent stationary increments, we obtain similarly as in \reff{martin},  for all $0\leq t\leq T$, 
\beqs
\E\big[  M_T(x,\theta^*,\alpha^*) \big|  \Fc_t \big] &=&  M_t(x,\theta^*,\alpha^*) \E[ M_{T-t}(0,\theta^*,\alpha^*)]  \\
&=&  M_t(x,\theta,a), 
\enqs
which shows that the process $M(x,\theta^*,\alpha^*)$  is a $(\P,\F)$-martingale.  Therefore, from the optional sampling theorem at the (bounded) 
stopping time $\tau_{x,\alpha^*}\wedge T$, we have 
\beq
e^{-\theta^* x}  &=& M_0(x,\theta^*,\alpha^*) \; = \; \E\big[ M_{\tau_{x,\alpha^*}\wedge T}(x,\theta^*,\alpha^*) \big] \nonumber \\
&=&  \E\big[ M_{\tau_{x,\alpha^*}}(x,\theta^*,\alpha^*) 1_{\tau_{x,\alpha^*} \leq T} \big] +  
\E\big[ M_{T}(x,\theta^*,\alpha^*) 1_{\tau_{x,\alpha^*} > T} \big] \nonumber \\
&\geq & \E\big[ M_{\tau_{x,\alpha^*}}(x,\theta^*,\alpha^*) 1_{\tau_{x,\alpha^*} \leq T} \big], \nonumber
\enq
since  the process $M$ is nonnegative. By the monotone convergence theorem, we then get by sending $T$ to infinity into the previous inequality
\beqs
e^{-\theta^* x} &\geq & \E\big[ M_{\tau_{x,\alpha^*}}(x,\theta^*,\alpha^*) 1_{\tau_{x,\alpha^*} < \infty} \big] \; = \; 
\E\big[ M_{\tau_{x,\alpha^*}}(x,\theta^*,\alpha^*) \big| \tau_{x,\alpha^*} <  \infty \big]  \P[\tau_{x,\alpha^*} < \infty],
\enqs
from Bayes formula. Thus,  we get
\beqs
\psi(x,\alpha^*) &=&   \P[\tau_{x,\alpha^*} < \infty]  
\; \leq \;  \frac{e^{-\theta^* x}}{\E\big[ M_{\tau_{x,\alpha^*}}(x,\theta^*,\alpha^*) \big| \tau_{x,\alpha^*} <  \infty \big]}. 
\enqs
Now, by definition of the time of ruin,  $V_{\tau_{x,\alpha^*}}^{x,\alpha^*}$ is nonpositive and so 
$M_{\tau_{x,\alpha^*}}(x,\theta^*,\alpha^*)$ $\geq$ $1$ a.s.  on $\{\tau_{x,\alpha^*}<\infty\}$.  We deduce that 
\beq \label{upperpsi}
\psi^*(x) \; \leq \; \psi(x,\alpha^*) & \leq & e^{-\theta^*x}. 
\enq

In order to state a lower bound on the minimal ruin probability, we proceed as follows. We apply It\^o's formula to the process 
$M(x,\theta^*,\alpha)$ for arbitrary $\alpha$ $\in$ $\Ac$~: 
\beqs
\frac{dM_t(x,\theta^*,\alpha)}{M_{t^-}(x,\theta^*,\alpha)} &=& \big( -\theta^*(p+b\alpha_t) + \frac{1}{2}(\theta^*)^2\sigma^2\alpha_t^2)dt 
- \theta^*\sigma dW_t + (e^{\theta^* Y_{N_t}}-1)dN_t. 
\enqs
Observing that $\gamma_{_Y}(\theta)$ $=$ $\E[e^{\theta Y_{N_t}}-1]$, we rewrite as
\beqs
\frac{dM_t(x,\theta^*,\alpha)}{M_{t^-}(x,\theta^*,\alpha)} &=& \big( -\theta^*(p+b\alpha_t) + \frac{1}{2}(\theta^*)^2\sigma^2\alpha_t^2 
+ \lambda \gamma_{_Y}(\theta^*))dt  - \theta^*\sigma dW_t  \\
& & \;\;\; + (e^{\theta^* Y_{N_t}}-1)dN_t - \lambda \E[e^{\theta^* Y_{N_t}}-1]dt \\
&=& f(\theta^*,\alpha_t) dt  - \theta^*\sigma dW_t + d \tilde N_t,
\enqs
where $\tilde N_t$ $=$ $\int_0^t (e^{\theta^* Y_{N_u}}-1)dN_u$ $-$ $\int_0^t \lambda \E[e^{\theta^* Y_{N_u}}-1]du$.  By using the martingale property of 
$N_t-\lambda t$, we can check that $\tilde N$ is a martingale. Since  $f(\theta^*,\alpha_t)$ $\geq$ $0$ a.s. for all $t$ (see \reff{fpos}), we deduce that 
$M(x,\theta^*,\alpha)$ is a (local) submartingale.  To go to a true submartingale, we need some additional assumption on the distribution of the claim size. 
Actually, we can prove that under the following uniform exponential tail distribution
\beq \label{uniexp}
\sup_{z\geq 0} \E[e^{-\theta^*(z-Y_1)} | Y_1 >  z] \; = \; \sup_{z\geq 0} \frac{\int_z^\infty e^{-\theta^*(z-y)}dF_Y(y)} {\int_z^\infty dF_Y} 
& < & \infty,
\enq
the process  $M(x,\theta^*,\alpha)$ is an uniformly integrable submartingale. Therefore, from the optional sampling theorem at the (bounded) 
stopping time $\tau_{x,\alpha}\wedge T$, we have 
\beq
e^{-\theta^* x}  \; = \;  M_0(x,\theta^*,\alpha) & \leq & \E\big[ M_{\tau_{x,\alpha}\wedge T}(x,\theta^*,\alpha) \big] \nonumber \\
&=&  \E\big[ M_{\tau_{x,\alpha}}(x,\theta^*,\alpha) 1_{\tau_{x,\alpha} \leq T} \big] +  
\E\big[ M_{T}(x,\theta^*,\alpha) 1_{\tau_{x,\alpha} > T} \big].  \label{intersub}
\enq
We now claim that $M_T(x,\theta^*,\alpha)$ converges a.s. on $\{\tau_{x,\alpha}=\infty\}$ to zero as $T$ goes to infinity. First, we know from Doob's submartingale convergence theorem that  $\lim_{T\rightarrow\infty}M_T(x,\theta^*,\alpha)$ exists a.s.,  hence also 
$\lim_{T\rightarrow\infty}V_T(x,\alpha)$. Since the expectation of a jump size $E[Y_1]$ is positive, there exists $y$ $>$ $0$ s.t. $\P[Y_1>y]$ $>$ $0$.  By independence of the jump sizes in the compound Poisson process of the risk reserve, it is then an easy exercice to see that with probability $1$, a jump of size greater than $y$ occurs infinitely often on $[0,\infty)$. On the other hand, the stochastic integral due to the invesment strategy $\alpha$ in the stock price, is a.s. continuous, and so cannot compensate the jumps of the compound Poisson process, greater than $y$, which will occur infinitely often a.s. It follows that on $\{\tau_{x,\alpha}=\infty\}$ 
(where ruin does not occur), $V_T^{x,\alpha}$  cannot converge to a finite value with positive probability. Therefore on 
$\{\tau_{x,\alpha}=\infty\}$, we have $\lim_{T\rightarrow\infty}V_T(x,\alpha)$  $=$ $\infty$ and thus, since $\theta^*$ $>$ $0$,  
$\lim_{T\rightarrow\infty}M_T(x,\theta^*,\alpha)$ $=$ $0$ a.s.  As $T$ $\rightarrow$ $\infty$, we have then from the dominated convergence theorem ($M_T(x,\theta^*,\alpha)$ $\leq$ $1$ on $\{\tau_{x,\alpha} > T\}$) for the second term in \reff{intersub}, and 
by the monotone convergence theorem for the first term,
\beqs
e^{-\theta^* x} & \leq &  \E\big[ M_{\tau_{x,\alpha}}(x,\theta^*,\alpha) 1_{\tau_{x,\alpha} < \infty} \big] 
 \; = \; 
\E\big[ M_{\tau_{x,\alpha}}(x,\theta^*,\alpha) \big| \tau_{x,\alpha} <  \infty \big]  \P[\tau_{x,\alpha} < \infty],
\enqs
and so 
\beq \label{phialpha}
\psi(x,\alpha) &=&   \P[\tau_{x,\alpha} < \infty]  
\; \geq \;  \frac{e^{-\theta^* x}}{\E\big[ M_{\tau_{x,\alpha}}(x,\theta^*,\alpha) \big| \tau_{x,\alpha} <  \infty \big]}. 
\enq
We finally prove that $\E\big[ M_{\tau_{x,\alpha}}(x,\theta^*,\alpha) \big| \tau_{x,\alpha} <  \infty \big]$ is bounded by a constant independent of $\alpha$ $\in$ $\Ac$. Fix some arbitrary $\alpha$ $\in$ $\Ac$ and set for shorthand notation $\tau$ $=$ $\tau_{x,\alpha}$ 
the time of ruin of the wealth process $V^{x,\alpha}$.  First observe that ruin $\{\tau <\infty\}$  occurs either through the brownian motion, i.e. on $\{\tau <\infty, V_{\tau^-}^{x,\alpha}=0\}$, and in this case $V_{\tau}^{x,\alpha}$ $=$ $0$ and so $M_{\tau}(x,\theta^*,\alpha)$ $=$ $1$, or through a jump, i.e. on 
$\{\tau <\infty, V_{\tau^-}^{x,\alpha} > 0\}$,  
and in this case $V_{\tau}^{x,\alpha}$ $<$ $0$ and so $M_{\tau}(x,\theta^*,\alpha)$ $>$ $1$.  Hence, 
\beq
\E\big[ M_{\tau}(x,\theta^*,\alpha) \big| \tau <  \infty \big] & \leq & 
\E\big[ M_{\tau}(x,\theta^*,\alpha) \big| \tau  <  \infty, V_{\tau^-}^{x,\alpha} > 0 \big] \nonumber \\
&=&  \E\big[ e^{-\theta^*(V_{\tau^-}^{x,\alpha} - Y_{N_\tau})} \big| \tau  <  \infty, V_{\tau^-}^{x,\alpha} > 0 \big]. \label{vinter}
\enq 
Let $H^{x,\alpha}(dt,dz)$ denote the joint distribution of $\tau$ and $V_{\tau^-}^{x,\alpha}$ conditional on the event 
$\{\tau  <  \infty, V_{\tau^-}^{x,\alpha} > 0\}$ that ruin occurs  through a jump.  Given $\tau$ $=$ $t$ and $V_{\tau^-}^{x,\alpha}$ $=$ $z$ 
$>$ $0$, a claim $Y_{N_\tau}$ occurs at time $t$ and has distribution $dF_{_Y}(y)/\int_z^\infty dF_{_Y}$ for $y$ $>$ $z$. Hence
\beq
\E\big[ e^{-\theta^*(V_{\tau^-}^{x,\alpha} - Y_{N_\tau})} \big| \tau  <  \infty, V_{\tau^-}^{x,\alpha} > 0 \big] &=& 
\int_0^\infty \int_0^\infty H^{x,\alpha}(dt,dz) \int_z^\infty e^{-\theta^*(z-y)} \frac{dF_{_Y}(y)}{\int_z^\infty dF_{_Y}} \nonumber \\
& \leq & \sup_{z\geq 0}  \int_z^\infty e^{-\theta^*(z-y)} \frac{dF_{_Y}(y)}{\int_z^\infty dF_{_Y}}  \; < \; \infty, \label{minter}
\enq
by assumption \reff{uniexp}. By setting 
\beqs
C &=& \frac{1}{\sup_{z\geq 0}  \int_z^\infty e^{-\theta^*(z-y)} \frac{dF_{_Y}(y)}{\int_z^\infty dF_{_Y}} } \; = \; 
\inf_{z\geq 0} \frac{\int_z^\infty dF_{_Y}}{\int_z^\infty e^{-\theta^*(z-y)} dF_{_Y}(y) } \; \in \; (0,1],
\enqs
we then have from \reff{phialpha}-\reff{vinter}-\reff{minter}, for all $x$ $\geq$ $0$, 
\beqs
\psi(x,\alpha) & \geq & C e^{-\theta^* x}, \;\;\; \forall \alpha\in \Ac.
\enqs
Together with the upper bound \reff{upperpsi}, this completes the proof of Theorem \ref{theoinsurfin}.

\section{Large deviations and rare event simulation in option pri\-cing}

\setcounter{equation}{0} \setcounter{Assumption}{0}
\setcounter{Theorem}{0} \setcounter{Proposition}{0}
\setcounter{Corollary}{0} \setcounter{Lemma}{0}
\setcounter{Definition}{0} \setcounter{Remark}{0}

\subsection{Importance sampling and large deviations approximations}  \label{secISlar}

In this paragraph, we show how to use large deviations approximation via importance sampling for 
Monte-carlo computation of expectations arising  in option pricing. In the context of continuous-time models, we are  interested in the 
computation of 
\beqs
I_g &=& \E\Big[ g(S_t,0\leq t\leq T) \Big], 
\enqs
where $S$ is the underlying asset price, and $g$ is the payoff of the option, eventually path-dependent, i.e. depending on the 
path process $S_t$, $0\leq t\leq T$.  The Monte-Carlo approximation technique consists in simulating $N$ independent 
sample paths $(S_t^i)_{0\leq t\leq T}$, $i$ $=$ $1,\ldots,N$,  in the distribution of $(S_t)_{0\leq t\leq T}$, and approximating the 
required expectation by the sample mean estimator~: 
\beqs
I_g^N &=& \frac{1}{N} \sum_{i=1}^N g(S^i). 
\enqs
The consistency of this estimator is ensured by the law of large numbers, while the error approximation is given by the variance of this estimator from the central limit theorem~:  the lower is the variance of $g(S)$, the better is the approximation for a given number $N$ of simulations.  As already mentioned in the introduction,  the basic principle of importance sampling is to reduce variance by changing probability measure from which paths are generated. Here, the idea is to change the distribution of the  price process to be simulated in order to take into account the specificities of the payoff function $g$.   
We focus in this section in the importance sampling technique within the context of diffusion models, and then show how to obtain an optimal change of measure by a large deviation approximation of the required expectation.

\subsubsection{Importance sampling for diffusions via Girsanov's theorem} \label{paragir}

We briefly describe  the importance sampling variance reduction technique for diffusions. Let $X$ be a $d$-dimensional diffusion process 
governed by
\beq \label{diffX}
dX_s &=& b(X_s) ds  + \Sigma(X_s) dW_s,
\enq
where $(W_t)_{t\geq 0}$ is a $n$-dimensional brownian motion on a filtered probability space $(\Omega,\Fc,\F = (\Fc_t)_{t\geq 0},\P)$,  and the Borel functions $b$, $\Sigma$ satisfy the usual Lipschitz condition ensuring the existence of a strong solution to the s.d.e. \reff{diffX}.  We denote by $X_s^{t,x}$ 
the solution to \reff{diffX} starting fom $x$ at time $t$, and we define the function~: 
\beqs
v(t,x) &=& \E\Big[ g(X_s^{t,x},t\leq s \leq T)\Big], \;\;\; (t,x) \in [0,T]\times\R^d. 
\enqs 
Let $\phi$ $=$ $(\phi_t)_{0\leq t\leq T}$ be an $\R^d$-valued adapted process such that the process
\beqs
M_t &=& \exp\Big( - \int_0^t \phi_u' dW_u - \frac{1}{2} \int_0^t |\phi_u|^2 du \Big), \;\;\; 0 \leq t\leq T,  
\enqs
is a martingale, i.e. $\E[M_T]$ $=$ $1$. This is ensured for instance under the Novikov criterion~: 
$\E\big[\exp\big(\frac{1}{2}\int_0^T |\phi_u|^2du\big)\big]$ $<$ $\infty$.  In this case, one can define a probability measure $\Q$ equivalent to $\P$ on 
$(\Omega,\Fc_T)$ by~: 
\beqs
\frac{d\Q}{d\P} &=& M_T.
\enqs
Moreover, by Girsanov's theorem, the process $\hat W_t$ $=$ $W_t+\int_0^t \phi_u du$, $0\leq t\leq T$, is a brownian motion under $\Q$, and the dynamics of $X$ under $\Q$ is given by
\beq \label{XsousQ}
dX_s &=& \big( b(X_s) - \Sigma(X_s) \phi_s\big) ds +  \Sigma(X_s) d\hat W_s. 
\enq
From Bayes formula, the expectation of interest can be written as
\beq \label{vsousQ}
v(t,x) &=& \E^{\Q}\Big[ g(X_s^{t,x},t\leq s \leq T) L_T \Big],
\enq
where  $L$ is the $\Q$-martingale
\beq \label{defL}
L_t \; = \; \frac{1}{M_t} &=&  \exp\Big( \int_0^t \phi_u' d\hat W_u - \frac{1}{2} \int_0^t |\phi_u|^2 du \Big), \;\;\; 0 \leq t\leq T.   
\enq 
The expression \reff{vsousQ} suggests, for any choice of $\phi$,  an alternative Monte-Carlo estimator for $v(t,x)$ with
\beqs
I_{g,\phi}^N(t,x)  &=& \frac{1}{N} \sum_{i=1}^N g(X^{i,t,x}) L_T^i,
\enqs 
by simulating $N$  independent sample paths $(X^{i,t,x})$ and $L_T^i$ of $(X^{t,x})$ and $L_T$  under $\Q$ given by \reff{XsousQ}-\reff{defL}.   
Hence, the change of probability measure through the choice of $\phi$ leads to a modification of the drift process in the simulation of $X$. 
The variance reduction technique consists in determining a process $\phi$, which induces a smaller variance for the  corresponding estimator $I_{g,\phi}$ than the initial one $I_{g,0}$.  The two next paragraphs present two approaches leading to the construction of such  processes $\phi$. In the first approach developed in \cite{foulastou97}, the process $\phi$ is stochastic, and requires an approximation of the expectation of interest. 
In the second approach due to \cite{glaheihas99}, the process $\phi$ is deterministic and derived through a simple optimization problem.  Both approaches rely on asymptotic results from the theory of large deviations.

\subsubsection{Option pricing approximation with  a Freidlin-Wentzell large deviation principle} \label{paraISbar}

We are looking for a stochastic process $\phi$, which allows to reduce (possibly to zero!) the variance of the corresponding estimator. 
The heuristics for achieving this goal is based on the following argument. Suppose for the moment that the payoff 
$g$ depends only on the terminal value $X_T$. Then, by applying It\^o's formula to  the $\Q$-martingale $v(s,X_s^{t,x})L_s$ between $s$ $=$ $t$ and $s$ $=$ $T$, we obtain~: 
\beqs
g(X_T^{t,x})L_T &=&  v(t,x) L_t + \int_t^T L_s \big( D_x v(s,X_s^{t,x})'\Sigma(X_s^{t,x}) + v(x,X_s^{t,x}) \phi_s'\big) d\hat W_s. 
\enqs
Hence, the variance of  $I_{g,\phi}^N(t,x)$ is given by
\beqs
Var_{\Q}(I_{g,\phi}^N(t,x)) &=& \frac{1}{N}  \E^{\Q} \Big[ \int_t^T L_s^2 \big| D_x v(s,X_s^{t,x})'\Sigma(X_s^{t,x}) + v(x,X_s^{t,x}) \phi_s'\big|^2 ds\Big]. 
\enqs
The choice of the process $\phi$ is motivated by the following remark. If the function $v$ were known, then one could vanish the variance by choosing 
\beq \label{phiopt}
\phi_s \; = \; \phi_s^* &=& - \frac{1}{v(s,X_s^{t,x})} \Sigma'(X_s^{t,x}) D_x v(s,X_s^{t,x}), \;\; t \leq s \leq T.  
\enq
Of course, the function $v$ is unknown (this is precisely what we want to compute), but this suggests to use a process $\phi$  from the above formula with an approximation of the function $v$. We may then reasonably hope to reduce  the variance, and also to use such a method for more general payoff  functions, possibly path-dependent.  We  shall use a large deviations approximation for the function $v$.

The  basic idea  for the use of large deviations approximation to the expectation function $v$ is the following.   
Suppose the option of interest, characterized by its payoff function $g$,  has a  low probability of exercice, e.g. it is deeply out the money.  Then, a large proportion of simulated paths 
end up out of the exercice domain,  giving no contribution to the Monte-carlo estimator but  increasing  the variance.  In  order to reduce the variance, 
it is interesting  to change of drift in the simulation of  price process to make the domain exercice more likely.  This is achieved with a large deviations 
approximation of the process  of interest in the asymptotics of small diffusion term~: such a result is known in the literature 
as  Freidlin-Wentzell sample path large deviations  principle.  Equivalently, by  time-scaling, this amounts to large deviation approximation 
of the process in small time, studied by Varadhan.

To illustrate our purpose, let us consider the case of an up-in bond, i.e. an option that pays one unit of num\'eraire iff the underlying asset   reached 
a given  up-barrier $K$. Within a stochastic volatility model $X$ $=$ $(S,Y)$  as in \reff{diffX} and given by~: 
\beq
dS_t &=& \sigma(Y_t) S_tdW_t^1 \label{dynS} \\
dY_t &=& \eta(Y_t) dt  + \gamma(Y_t) dW_t^2, \;\;\; d<W_1,W_2>_t  \; = \; \rho dt, \label{dynY} 
\enq
its price is then given by
\beqs
v(t,x) &=& \E\big[ 1_{\max_{t\leq u\leq T} S_u^{t,x} \geq K} \big]  \; = \;  \P[ \tau_{_{t,x}} \leq T ], \;\; t \in [0,T], \; x = (s,y) \in (0,\infty)\times\R, 
\enqs 
where 
\beqs
\tau_{_{t,x}} &=& \inf \big\{ u\geq t~: X_u^{t,x} \notin \Gamma \big\}, \;\;\;\;\;  \Gamma \; = \; (0,K) \times\R. 
\enqs
The event  $\{ \max_{t\leq u\leq T}  S_u^{t,x} \geq K \}$ $=$ $\{\tau_{_{t,x}} \leq T\}$ is rare when  $x$ $=$ $(s,y)$ $\in$ $\Gamma$, i.e. $s$ $<$ $K$ (out the money option) and the time to maturity $T-t$ is small.  The large deviations asymptotics for the exit probability  $v(t,x)$ in small time to maturity  $T-t$ 
is  provided by  the Freidlin-Wentzell and Varadhan theories.  
Indeed, we see from the time-homogeneity of the coefficients  of the diffusion and by time-scaling that  we may write $v(t,x)$ $=$ $w_{_{T-t}}(0,x)$, where for $\eps$ $>$ $0$,   $w_\eps$ is the function defined on 
$[0,1]\times (0,\infty)\times\R$ by
\beqs
w_\eps(t,x)  &=&  \P[\tau_{_{t,x}}^\eps  \leq 1 ],  
\enqs 
and $X^{\eps,t,x}$   is the solution to 
\beqs
dX^\eps_s &=& \eps b(X_s^\eps) ds + \sqrt{\eps} \Sigma(X_s^\eps) dW_s, \;\;\; X_t^\eps \; = \; x. 
\enqs
and $\tau_{_{t,x}}^\eps$ $=$  $\inf \big\{ s \geq t~: X_s^{\eps,t,x} \notin \Gamma \big\}$.  From the large deviations result  \reff{larfrei} stated in  
paragraph \ref{paraprin}, we have~: 
\beqs
\lim_{t\nearrow T} - (T-t) \ln v(t,x)  &=& V_0(0,x),
\enqs
where 
\beqs
V_0(t,x) &=& \inf_{x(.) \in \Ac(t,x)}  \int_t^1   \frac{1}{2} \dot x(u)' M(x(u)) \dot x(u)  du,  \;\;\;  (t,x) \in [0,1)\times\Gamma, 
\enqs
$\Sigma(x)$ $=$ $\left( \begin{array}{cc} \sigma(x) & 0 \\ 0 &  \gamma(x) \end{array} \right)$ is the diffusion matrix of $X$ $=$ $(S,Y)$, 
$M(x)$ $=$ $(\Sigma\Sigma'(x))^{-1}$, and 
\beqs
\Ac(t,x) &=& \big\{ x(.) \in H^1_{loc}([0,1],(0,\infty)\times\R)~:   \;  x(t) = x \; \mbox{ and } \; \tau(x) \leq  1\big\}. 
\enqs
Here,  for an absolutely continuous  function $x(.)$ on $[0,1)$ and valued in $(0,\infty)\times\R$, we denote $\dot{x}(u)$  its time derivative, and 
$\tau(x)$ the exit time of $x(.)$ from $\Gamma$. 

We also have another interpretation of  the positive function $V_0$ in terms of Riemanian distance on $\R^d$ associated to the metric  
$M(x)$ $=$ $(\Sigma\Sigma'(x))^{-1}$.  By denoting $L_0(x)$ $=$ $\sqrt{2V_0(0,x)}$,  one can prove (see  \cite{laslio95}) that  $L_0$ is 
the unique viscosity solution to the eikonal equation
\beqs
(D_x L_0)' \Sigma\Sigma'(x) D_x L_0 &=& 1, \;\;\; x \in \Gamma \\
L_0(x) &=& 0, \;\;\; x \in \partial\Gamma
\enqs
and that it may be represented as  
\beq \label{repdis}
L_0(x)  &=& \inf_{z \in \partial\Gamma} L_0(x,z), \;\;\;  x \in \Gamma, 
\enq
where
\beqs
L_0(x,z) &=& \inf_{x(.) \in A(x,z)} \int_0^1   \sqrt{\dot{x}(u)' M(x(u)) \dot{x}(u)} du,
\enqs
and $A(x,z)$ $=$ $\big\{ x(.) \in H^1_{loc}([0,1],(0,\infty)\times\R)~:   \;  x(0) = x \; \mbox{ and } \;  x(1) = z \big\}$.  Hence, the function $L_0$ 
can be computed  either by 
the numerical  resolution of the eikonal equation or by using the representation \reff{repdis}.   $L_0(x)$ is interpreted as the minimal  
length (according to the metric $M$) of the path $x(.)$ allowing to  reach  the boundary $\partial\Gamma$ from $x$.  From the above 
large deviations result, which  is written as
\beqs
\ln v(t,x) & \simeq &    - \frac{L_0^2(x)}{2(T-t)}, \;\;\; \mbox{ as }  \; T-t \rightarrow 0,  
\enqs
and the expression \reff{phiopt} for the optimal theoretical $\phi^*$, we use a change of probability measure with 
\beqs
\phi(t,x) &=& \frac{L_0(x)}{T-t} \Sigma'(x) D_x L_0(x). 
\enqs 
Such a process $\phi$ may also appear interesting to use in a more general framework than up-in bond~: one can use it for computing 
any option whose exercice domain looks similar to the up and in bond.  We also expect that the variance reduction is more significant as the exercice probability is low, i.e. for deep out the money options.  In the particular case of the Black-Scholes model, i.e. $\sigma(x)$ $\sigma s$, we have 
\beqs
L_0(x) &=& \frac{1}{\sigma} \big| \ln\big(\frac{s}{K}\big) \big|,
\enqs
and so 
\beqs
\phi(t,x) &=& \frac{1}{\sigma(T-t)}   \ln(\frac{s}{K}\big), \;\;\; s < K.  
\enqs

\subsubsection{Change of drift via Varadhan-Laplace principle}

We describe here a method due to \cite{glaheihas99}, which, in contrast with the above approach,  does not require the knowledge of the option price. This method restricts to deterministic changes of drift  over discrete time steps.  Hence, the diffusion model $X$ for 
state variables (stock price, volatility) 
is simulated (eventually using an Euler scheme if needed) on a discrete time grid $0=t_0$ $<$ $t_1$ $<\ldots<$ $t_m$ $=$ $T$~: the increment of the brownian motion from $t_{i-1}$ to $t_i$ is simulated as $\sqrt{t_i-t_{i-1}}Z_i$, where 
$Z_1,\ldots,Z_m$  are i.i.d. $n$-dimensional standard normal random vectors.  We denote by $Z$ the concatenation of the $Z_i$ 
into a single vector of lengh $l$ $=$ $mn$.  Each outcome of $Z$ determines a path of state variables. Let $G(Z)$ denote the payoff derived from $Z$, and our aim is to compute the (path-dependent) option price $\E[G(Z)]$. For example, in the case of the Black-Scholes model for the stock price $S$, we have
\beqs
S_{t_i} &=& S_{t_{i-1}} \exp\Big(-\frac{\sigma^2}{2}(t_i-t_{i-1}) + \sigma \sqrt{t_i-t_{i-1}}Z_i\Big),
\enqs
and the payoff of the Asian option is 
\beqs
G(Z) &=& G(Z_1,\ldots,Z_m) \; = \; \Big( \frac{1}{m} \sum_{i=1}^m S_{t_i} - K\Big)_+
\enqs

We apply importance sampling by changing the mean of $Z$ from $0$ to some vector $\mu$ $=$ $(\mu_1,\ldots,\mu_m)$. We denote 
$\P_\mu$ and $\E_\mu$ the probability and expectation when $Z$ $\sim$ $\Nc(\mu,I_m)$.  Notice that with the notations of paragraph \ref{paragir},  this corresponds to a piecewise constant process $\phi$ s.t. $\phi_t$ $=$ $-\mu_i/\sqrt{t_i-t_{i-1}}$ on $[t_{i-1},t_i)$. 
By Girsanov's theorem or here more simply from the likelihood ratio for normal random vectors, the corresponding unbiased estimator is then obtained by taking the average of independent replications of 
\beqs
\vartheta_\mu &=& G(Z) e^{-\mu'Z + \frac{1}{2}\mu'\mu},
\enqs
where $Z$ is sampled from $\Nc(\mu,I_m)$. We call $\vartheta_\mu$ a $\mu$-IS estimator.  In order to minimize over $\mu$ the variance of this estimator, it suffices to minimize its second moment, which is given by~: 
\beqs
M^2(\mu) &=& \E_\mu\Big[ G(Z)^2 e^{-2\mu'Z +  \mu'\mu} \Big] \; = \; \E\Big[ G(Z)^2 e^{-\mu'Z +  \frac{1}{2} \mu'\mu} \Big]
\enqs
We are then looking for an optimal $\mu$ solution to 
\beq \label{muopt1}
\inf_{\mu} M^2(\mu) \; : = \;  \inf_{\mu} \E\Big[ G(Z)^2 e^{-\mu'Z +  \frac{1}{2} \mu'\mu} \Big]. 
\enq
This minimization problem is, in general, a well-posed problem. Indeed, it is shown in \cite{aro04} that if $\P[G(Z)>0]$ $>$ $0$, and $\E[G(Z)^{2+\delta}]$ 
$<$ $\infty$ for some $\delta$ $>$ $0$, then $M^2(.)$ is a strictly convex function, and thus $\mu^*$ solution to \reff{muopt1} exists and is unique. 
This $\mu^*$ can be computed by solving (numerically) $\nabla M^2(\mu)$ $=$ $0$. This is the method adopted in \cite{aro04} with a 
Robbins-Monro stochastic algorithm.  We present here an approximate resolution of \reff{muopt1}  by means of large deviations approximation.  
For this,  assume that $G$ takes only nonnegative values, so that it is written as $G(z)$ $=$ $\exp(F(z))$, with the convention that $F(z)$ $=$ $-\infty$ if $G(z)$ $=$ $0$, and 
let us  consider the more general estimation problem where $Z$ is replaced by $Z^\eps$ $=$ $\sqrt{\eps} Z$ and we simultaneously scale the payoff by raising it to the power of $1/\eps$~: 
\beqs
v_\eps &=& \E[e^{\frac{1}{\eps}F(Z^\eps)}]
\enqs
The quantity of interest $\E[G(Z)]$ $=$ $\E[e^{F(Z)}]$ is $v_\eps$ for $\eps$ $=$ $1$. We embed the problem of estimating 
$v_1$ in the more general problem of estimating $v_\eps$ and  analyze the second moment of corresponding IS  estimators as 
$\eps$ is small, by means of Varadhan-Laplace principle. For any $\mu$, we consider $\mu_\eps$-IS estimator of $v_\eps$ with $\mu^\eps$ $=$ $\mu/\sqrt{\eps}$~: 
\beqs
\vartheta^\eps_{\mu} &=&  e^{\frac{1}{\eps}F(\sqrt{\eps}Z)} e^{-\mu_\eps'Z + \frac{1}{2}\mu_\eps'\mu_\eps} \; = \; 
e^{\frac{1}{\eps}(F(Z^\eps)-\mu'Z^\eps + \frac{1}{2}\mu'\mu)}
\enqs
where $Z$ is sampled from $\Nc(\mu_\eps,I_m)$ $=$ $\Nc(\mu/\sqrt{\eps},I_m)$. Its second moment is
\beqs
M_\eps^2(\mu) &=& \E_{\mu_\eps}\Big[  e^{\frac{1}{\eps}(2F(Z^\eps)-2\mu'Z^\eps + \mu'\mu)}  \Big] \; = \; 
\E\Big[  e^{\frac{1}{\eps}(2F(Z^\eps)-\mu'Z^\eps + \frac{1}{2} \mu'\mu})  \Big]
\enqs
Now, from Cramer's theorem, $(Z^\eps)_\eps$ satisfies a LDP with rate function $I(z)$ $=$ $\frac{1}{2}z'z$.  
Hence, under the condition that $F(z)$ $\leq$ $c_1+c_2z'z$ for some $c_2$ $<$ $1/4$, one can apply Varadhan's integral principle (see Theorem \ref{theovaradh}) to the function $z$ $\rightarrow$ $2F(z)-\mu'z + \frac{1}{2} \mu'\mu$, and get
\beq \label{optm2}
\lim_{\eps\rightarrow 0} \eps \ln M_\eps^2(\mu) &=& \sup_{z}\big[2F(z)-\mu'z + \frac{1}{2} \mu'\mu -  \frac{1}{2}z'z\big]. 
\enq
This  suggests  to  search for a $\mu$ solution to the problem~:   
\beq \label{minmax}
\inf_\mu \sup_{z}\big[2F(z)-\mu'z + \frac{1}{2} \mu'\mu -  \frac{1}{2}z'z\big]. 
\enq
This min-max problem may be reduced to a simpler one. Indeed, assuming that the conditions of the min/max theorem hold, then the $\inf$ and $\sup$ 
can be permuted, and we find
\beq
\inf_\mu \sup_{z}\big[2F(z)-\mu'z + \frac{1}{2} \mu'\mu -  \frac{1}{2}z'z\big] &=&  
\sup_{z} \Big[ \inf_\mu\big( -\mu'z + \frac{1}{2} \mu'\mu \big) + 2F(z)  -  \frac{1}{2}z'z \Big] \nonumber \\
&=& 2 \sup_z [F(z) -  \frac{1}{2}z'z \big]. \label{eq2opt} 
\enq
Actually, under suitable convexity conditions on $F$ and its domain, one can show (see \cite{glaheihas99}) that \reff{eq2opt} holds.
Furthermore,  if $\hat z$ is solution to  
\beq \label{maxmu}
\sup_z \big(F(z) - \frac{1}{2}z'z\big),
\enq
then a solution $\hat\mu$ solution to  \reff{minmax} should be identified with the conjugate of $\hat z$, via 
$\inf_\mu\big( -\mu'\hat z + \frac{1}{2} \mu'\mu \big)$ that is $\hat\mu$ $=$ $\hat z$.   
The solution to problem \reff{maxmu} has also the following interpretation.  From heuristic arguments of importance sampling (see the introduction), an optimal effective importance sampling density should assign high probability to regions on which the product of the integrand payoff and the original 
density is large. For our problem,  this product is proportional to 
\beqs
e^{F(z)-\frac{1}{2}z'z},
\enqs
since $\exp(-z'z/2)$ is proportional to the standard normal density. This suggests to choose $\mu$ $=$ $\hat\mu$  solution to \reff{maxmu}. Another heuristics indication for  the choice of \reff{maxmu} is based on the following argument. Assume that $F$ is $C^1$ on its domain, and the maximum 
$\hat\mu$ in  \reff{maxmu} is attained in the interior of the domain, so that  it solves the fixed point equation~:
\beq \label{fixpoint}
\nabla F(\hat\mu) &=& \hat\mu. 
\enq
By using a first-order Taylor approximation of $F$ around the mean $\hat\mu$ of $Z$ under $\P_{\hat\mu}$, we may approximate the estimator as
\beq \label{Taylor}
\vartheta_{\hat\mu} & = & e^{F(Z)-\hat\mu'Z+\frac{1}{2}\hat\mu'\hat\mu}  \; \simeq \; 
e^{F(\hat\mu) + \nabla F(\hat\mu)'(Z-\hat\mu) -\hat\mu'Z+\frac{1}{2}\hat\mu'\hat\mu }
\enq
Hence, for the choice of $\hat\mu$ satisfying \reff{fixpoint}, the expression of the r.h.s. of \reff{Taylor} collapses to a constant with no dependence on $Z$. Thus, applying importance sampling with such a $\hat\mu$ would produce a zero-variance estimator if \reff{Taylor} holds exactly, e.g. if $F$ is linear, and it should produce a low-variance estimator if \reff{Taylor} holds only approximately.

The choice of $\mu$ $=$ $\hat\mu$ solution to \reff{maxmu} leads also to an asymptotically optimal IS-estimator in the following sense. 
First, notice  that for any $\mu$, we have by Cauchy-Schwarz's inequality~: $M_\eps^2(\mu)$ $\geq$ $(v_\eps)^2$, and so
\beqs
\lim_{\eps\rightarrow 0} \eps \ln M_\eps^2(\mu) & \geq & 2 \lim_{\eps\rightarrow 0} \eps \ln v_\eps 
\; = \; 2  \lim_{\eps\rightarrow 0} \eps \ln \E[e^{\frac{1}{\eps}F(Z^\eps)}]  
\enqs
From Varadhan's integral principle  applied to the function $z$ $\rightarrow$ $F(z)$, we thus deduce for any $\mu$,
\beqs
\lim_{\eps\rightarrow 0} \eps \ln M_\eps^2(\mu) & \geq & 2 \sup_{z}[F(z) -  \frac{1}{2}z'z\big] \; = \; 
2[F(\hat\mu) - \frac{1}{2}\hat\mu'\hat\mu].
\enqs
Hence, $2[F(\hat\mu) - \frac{1}{2}\hat\mu'\hat\mu]$ is the best-possible exponential decay rate for a $\mu_\eps$-IS estimator $\vartheta_\mu^\eps$.   
Now, by choosing $\mu$ $=$ $\hat\mu$, and from \reff{optm2}, \reff{eq2opt}, we have 
\beqs
\lim_{\eps\rightarrow 0} \eps \ln M_\eps^2(\hat\mu) &=& 2[F(\hat\mu) - \frac{1}{2}\hat\mu'\hat\mu],
\enqs 
 which shows that the $\hat\mu_\eps$-IS estimator $\vartheta_{\hat\mu}^\eps$ is asymptotically optimal. 
 
\begin{Remark}
{\rm From  the first-order approximation of $F$ in \reff{Taylor}, we see that  in order to obtain further variance reduction, it is natural to address the quadratic component of $F$. This can be achieved by a method of stratified sampling as developed  in \cite{glaheihas99}. 
}
\end{Remark}

\vspace{2mm}

Recently,  the above approach of \cite{glaheihas99} was extended in \cite{guarob06} to a continuous-time setting, where the optimal 
deterministic drift in the Black-Scholes model is characterized as the solution to a classical one-dimensional variational problem.

\subsection{Computation of barrier crossing probabilities and sharp large deviations}

In this paragraph, we present a simulation procedure for computing the probability that a diffusion process crosses pre-specified barriers in a given time interval $[0,T]$. Let $(X_t)_{t\in [0,T]}$ be a diffusion process in $\R^d$, 
\beqs
dX_t &=& b(X_t) dt + \sigma(X_t) dW_t
\enqs
and $\tau$ is the exit time of $X$ from some domain $\Gamma$ of $\R^d$, eventually depending on time~:   
\beqs
\tau &=& \inf\{ t \in [0,T]~: X_t \notin \Gamma(t) \},
\enqs  
with the usual convention that $\inf\emptyset$ $=$ $\infty$. Such a quantity appears typically in finance in the computation of barrier options, for 
example with a knock-out  option~:  
\beq \label{prixbar}
C_0 &=& \E\big[e^{-rT} g(X_T) 1 _{\tau > T}\big],
\enq  
with  $\Gamma(t)$ $=$ $(-\infty,U(t))$ in the case of single barrier options, and  $\Gamma(t)$ $=$ $(L(t),U(t))$, for double barrier options. Here, 
$L$, $U$ are real functions  $:$  $[0,\infty)$ $\rightarrow$ $(0,\infty)$ s.t. $L$ $<$ $U$.

The direct naive approach would consist  first of simulating the process $X$ on $[0,T]$ through a discrete Euler scheme of step size 
$\eps$ $=$ $T/n$ $=$  $t_{i+1}-t_i$, $i$ $=$ $0,\ldots,n$~: 
\beqs
\bar X_{t_i+1}^\eps &=& \bar X_{t_i}^\eps +  b(\bar X_{t_i}^\eps) \eps + \sigma(\bar X_{t_i}^\eps) (W_{t_{i+1}}- W_{t_i}),
\enqs 
and the exit time $\tau$ is approximated  by the first time the discretized process reaches the barrier~:
\beqs
\bar\tau^\eps & = & \inf\big\{ t_i~:  \bar X_{t_i}^\eps \notin \Gamma(t_i) \big\}.
\enqs
Then, the barrier option price $C_0$ in  \reff{prixbar} is approximated by  Monte-Carlo simulations of the quantity
\beqs
\bar C_0^\eps &=& \E\big[e^{-rT}g(\bar X_T^\eps) 1 _{\bar\tau^\eps > T}\big]. 
\enqs 
In this procedure,  one considers that the price diffusion is killed if there exists a value $\bar X_{t_i}^\eps$, which is out of the domain $\Gamma(t_i)$.  
Hence,  such an approach is suboptimal since it does not control the diffusion path between two successive dates $t_i$ and $t_{i+1}$~:  the     
diffusion path could have crossed the barriers and come back to the domain without being detected. In this case, one over-estimates 
the exit time probability of $\{\tau > T\}$. This suboptimality is confirmed by the property 
that the error between $C_0$ and $\bar C_0^\eps$ is of order $\sqrt{\eps}$, as shown in \cite{gob00},  instead of the usual order $\eps$ obtained for standard vanilla options.

In order to improve  the above  procedure, we need to determine the probability that the process $X$ crosses the barrier between discrete simulation times.  We then consider  the continuous Euler scheme
\beqs
\bar X_{t}^\eps &=& \bar X_{t_i}^\eps +  b(\bar X_{t_i}^\eps)(t-t_i)  + \sigma(\bar X_{t_i}^\eps) (W_t - W_{t_i}), \;\;\; t_i \leq t \leq t_{i+1},
\enqs
which evolves as a Brownian with drift between two time discretizations $t_i$,  $t_{i+1}$ $=$ $t_i+\eps$.  
Given a simulation path of $(\bar X_{t_i}^\eps)_i$,  and values $\bar X_{t_i}^\eps$ $=$ $x_i$,  $\bar X_{t_{i+1}}^\eps$ $=$ $x_{i+1}$,  we denote 
\beqs
p_i^\eps(x_i,x_{i+1}) &=& \P\Big[ \exists t \in [t_i,t_{i+1}]~:  \bar X_t^\eps \notin \Gamma(t_i) \big| ( \bar X_{t_i}^\eps,\bar X_{t_{i+1}}^\eps)=(x_i,x_{i+1})  \Big],  
\enqs
the exit probability of the Euler scheme conditionally on the simulated path values.  The correction Monte-Carlo procedure works then as follows~: 
with probability $p_i^\eps$ $=$ $p_i^\eps(\bar X_{t_i}^\eps,\bar X_{t_{i+1}}^\eps)$, we stop the simulation by considering that the diffusion is killed, and we set $\tau^\eps$ $=$ $t_i$;  with probability $1-p_i^\eps$, we carry on the simulation.  The approximation of \reff{prixbar} is thus computed by Monte-Carlo simulations of 
\beqs
C_0^\eps &=& \E\big[e^{-rT}(\bar X_T^\eps -K)_+ 1 _{\tau^\eps > T}\big].  
\enqs 
We then recover a rate of convergence  of order $\eps$ for  $C_0^\eps - C_0$, see \cite{gob00}.

The effective implementation of this corrected procedure requires the calculation of $p_i^\eps$. 
Notice that on the interval $[t_i,t_{i+1}]$,  the 
diffusion $\bar X^\eps$ conditionned to $\bar X_{t_i}^\eps$ $=$ $x_i$,  $\bar X_{t_{i+1}}^\eps$ $=$ $x_{i+1}$, is a  brownian bridge~:  it coincides in distribution with  the process
\beqs
\tilde B_t^{i,\eps} &=& x_i + \frac{t}{\eps}\big(x_{i+1} - x_i\big) + \sigma(x_i)\big(W_t -  \frac{t}{\eps} W_\eps), \;\;\;  0 \leq t \leq \eps,
\enqs
and so by time change $t$ $\rightarrow$ $t/\eps$, with the process
\beqs
Y_t^{i,\eps} \; := \; \tilde B_{\eps t}^{i,\eps} &=& x_i +  t  \big(x_{i+1} - x_i\big) 
+ \sqrt{\eps} \sigma(x_i)\big(W_t -  t  W_1), \;\;\;  0 \leq t \leq 1. 
\enqs
It is known that the process $Y^{i,\eps}$ is solution to the s.d.e.
\beqs
dY_t^{i,\eps} &=& - \frac{Y_t^{i,\eps} - x_{i+1}}{1-t} dt  + \sqrt{\eps} \sigma(x_i) dW_t, \;\;\; 0 \leq t < 1, \\
Y_0^{i,\eps} &=& x_i.  
\enqs
The probability $p_i^\eps$ can then be expressed as 
\beq \label{pifrei}
p_i^\eps(x_i,x_{i+1}) &=& \P[ \tau^{i,\eps} \leq 1 ],  \;\;\; \mbox{ where } \;\; 
\tau^{i,\eps} \; = \; \inf\big\{ t \geq 0~: Y_t^{i,\eps} \notin \Gamma(t_i+ \eps t) \big\}. 
\enq
In the case of a half-space, i.e. single constant barrier, one has an explicit expression for the exit probability of a Brownian bridge. 
For example, if $\Gamma(t)$ $=$  $(-\infty,U)$, we have 
\beqs
p_i^\eps(x_i,x_{i+1}) &=& \exp\Big( - \frac{I_{_U}(x_i,x_{i+1})}{\eps}  \Big), \;\; \mbox{ with } \;\; 
I_{_U}(x_i,x_{i+1}) \; = \; \frac{2}{\sigma^2(x_i)}(U-x_i)(U-x_{i+1}).  
\enqs 
In the general case, we do not have analytical expressions for $p_i^\eps$, and one has to rely on simulation techniques or asymptotic approximations. 
We shall here consider asymptotic techniques based on large deviations  and Freidlin-Wentzell theory.  
Let us illustrate this point  in the case of two time-dependent barriers, i.e. $\Gamma(t)$ $=$ $(L(t),U(t))$ for smooth barriers functions $L$ $<$ $U$. 
Problem \reff{pifrei} does not exactly fit into the Freidlin-Wentzell framework considered in paragraph \ref{paraprin}, but was adapted for Brownian 
bridges with time-dependent barriers in \cite{balcariov99}.  We then have the large deviation estimate for $p_i^\eps$~: 
\beqs
\lim_{\eps \rightarrow 0} \eps \ln p_i^\eps(x_i,x_{i+1}) &=& - I_{_{L,U}}(x_i,x_{i+1}),
\enqs
where $I_{_{L,U}}(x_i,x_{i+1})$ is the infimum of the functional 
\beqs
y(.)  & \longrightarrow & \frac{1}{2\sigma(x_i)^2 } \int_0^1 \big| \dot y(t)  + \frac{y(t)- x_{i+1}}{1-t} \big|^2 dt,
\enqs
over all absolutely continuous paths $y(.)$ on $[0,1]$ s.t. $y(0)$ $=$ $x_i$, and there exists some $t$ $\in$ $[0,1]$ for which 
$y(t)$ $\leq$ $L(t_i)$ or $y(t)$ $\geq$ $U(t_i)$.  This infimum is a classical problem of calculus of variations, which is explicitly solved and gives for any 
$x_i,x_{i+1}$ $\in$ $(L(t_i),U(t_i))$ (otherwise $I_{_{L,U}}(x_i,x_{i+1})$ $=$ $0$)~: 
\beqs
I_{_{L,U}}(x_i,x_{i+1}) &=& \left\{ \begin{array}{cc}  
                                                      \frac{2}{\sigma^2(x_i)}(U(t_i)-x_i)(U(t_i)-x_{i+1}) & \mbox{ if } \; x_i + x_{i+1} > L(t_i) + U(t_i) \\
                                                      \frac{2}{\sigma^2(x_i)}(x_i-L(t_i))(x_{i+1}-L(t_i)) & \mbox{ if } \; x_i + x_{i+1} < L(t_i) + U(t_i).   
                                                      \end{array}
                                                      \right.
\enqs
In order to remove the log estimate on $p_i^\eps$, we need  a sharper large deviation estimate, and this is analyzed by the results of \cite{flejam92} recalled in paragraph \ref{paraprin}.   More precisely, we have  
\beqs
p_i^\eps(x_i,x_{i+1}) &=& \exp\Big( - \frac{I_{_{L,U}}(x_i,x_{i+1})}{\eps} -  w_{_{L,U}}(x_i,x_{i+1}) \Big) (1+ O(\eps)),
\enqs 
where $w_{_{L,U}}(x_i,x_{i+1})$ is explicited in \cite{balcariov99} as
\beqs
w_{_{L,U}}(x_i,x_{i+1}) &=& \left\{ \begin{array}{cc}  
                                                      \frac{2}{\sigma^2(x_i)}(U(t_i)-x_i)U'(t_i) & \mbox{ if } \; x_i + x_{i+1} > L(t_i) + U(t_i) \\
                                                      \frac{2}{\sigma^2(x_i)}(x_i-L(t_i)) L'(t_i) & \mbox{ if } \; x_i + x_{i+1} < L(t_i) + U(t_i).   
                                                      \end{array}
                                                      \right.
\enqs
Some recent extensions of this large deviations approach to the computation of exit probabilities for multivariate Brownian bridge are studied in 
\cite{huhkol06}, which also gives applications for the estimation of default probabilities in credit risk models, and the pricing of credit default swaps.

\section{Large deviations in risk management}

\setcounter{equation}{0} \setcounter{Assumption}{0}
\setcounter{Theorem}{0} \setcounter{Proposition}{0}
\setcounter{Corollary}{0} \setcounter{Lemma}{0}
\setcounter{Definition}{0} \setcounter{Remark}{0}

\subsection{Large portfolio losses in  credit risk}

\subsubsection{Portfolio credit risk in a single factor normal copula model}

A basic problem in measuring portfolio credit risk is determining the distribution of losses from default over a fixed horizon. 
Credit portfolios are often large, including exposure to thousands of obligors, and the default probabilities of high-quality credits are extremely small. 
These features in credit risk context lead to consider rare but significant large loss events, and emphasis is put on  the small probabilities of large losses, as these are relevant for  calculation of value at risk and related risk measures.  

We use the following notation~:

\vspace{2mm}

$n$ $=$ number of obligors to which portfolio is exposed,

$Y_k$ $=$  default indicator ($=$ $1$ if default, $0$ otherwise) for $k$-th obligor,

$p_k$ $=$ marginal probability that $k$-th obligor defaults, i.e. $p_k$ $=$ $\P[Y_k=1]$,

$c_k$ $=$ loss resulting from default of the $k$-th obligor,

$L_n$ $=$ $c_1Y_1+\ldots+c_nY_n$ $=$ total loss from defaults.

\vspace{2mm}

\noindent We are interested in estimating tail probabilities $\P[L_n > \ell_n]$ in the limiting regime at increasingly high loss  thresholds $\ell_n$, 
and rarity of large losses  resulting from a  large number $n$ of obligors and multiple defaults.

For simplicity, we  consider a homogeneous portfolio where all $p_k$ are equal to $p$, and all $c_k$  are equal constant  to $1$.  
An essential feature for credit risk management is the mechanism used to model the dependence across sources of credit risk.    
The dependence among obligors is modelled by the dependence among the default indicators $Y_k$. This dependence is introduced through a normal copula model as follows~: each default indicator is represented as
\beqs
Y_k &=& 1_{\{X_k > x_k\}}, \;\;\; k = 1,\ldots,n,
\enqs
where $(X_1,\ldots,X_n)$ is a multivariate normal vector. Without loss of generality, we take each $X_k$ to have a standard normal distribution, and we 
choose $x_k$  to match the marginal default probability $p_k$, i.e. $x_k$ $=$ $\Phi^{-1}(1-p_k)$ $=$ $-\Phi^{-1}(p_k)$, with $\Phi$ cumulative normal distribution.  We also denote $\varphi$ $=$ $\Phi'$ the density of the normal distribution. 
The correlations along the $X_k$, which determine the dependence among the $Y_k$, are specified through a single factor model of the form~: 
\beq \label{single}
X_k &=& \rho Z + \sqrt{1-\rho^2} \eps_k,  \;\;\; k =1,\ldots,n. 
\enq
where $Z$ has the  standard normal distribution $\Nc(0,1)$, $\eps_k$ are independent $\Nc(0,1)$ distribution, and $Z$ is independent of  $\eps_k$,  $k$ $=$ $1,\ldots,n$.  $Z$ is called systematic risk factor (industry, regional risk factors for example ...), and $\eps_k$ is an idiosyncratic risk associated with the $k$-th obligor. The constant $\rho$ in $[0,1)$ is a factor loading on the single factor $Z$, and assumed here to be identical for all obligors.  
We shall distinguish the case of independent obligors ($\rho$ $=$ $0$), and dependent obligors $(\rho$ $>$ $0$). 
More general multivariate factor models with inhomogeneous obligors are studied in \cite{glakansha06}.

\subsubsection{Independent obligors}

In this case, $\rho$ $=$ $0$, the default indicators $Y_k$  are i.i.d. with Bernoulli distribution of parameter $p$, and $L_n$ is a binomial distribution of parameters $n$ and $p$.  By the law of large numbers, $L_n/n$ converges to $p$.  Hence, in order that the loss event  $\{L_n \geq l_n\}$ becomes rare (without being trivially impossible), we let $l_n/n$ approach $q$ $\in$ $(p,1)$. It is then appropriate to specify $l_n$ $=$ $nq$ with $p< q < 1$. 
From Cramer's theorem and the expressions of the c.g.f. of the Bernoulli distribution and its Fenchel-Legendre transform, we obtain the large deviation 
result for the loss probability~: 
\beqs
\lim_{n\rightarrow} \frac{1}{n} \ln \P[ L_n \geq nq ] &=& - q \ln\big(\frac{q}{p}\big) - (1-q)  \ln\big(\frac{1-q}{1-p}\big) \; < \; 0. 
\enqs

\begin{Remark}
{\rm  By denoting $\Gamma(\theta)$ $=$ $\ln(1-p+pe^\theta)$ the c.g.f. of  $Y_k$, we have an IS (unbiased) estimator of  $\P[L_n\geq nq]$ by taking the average  of independent replications of 
\beqs
\exp(-\theta L_n + n\Gamma(\theta)) 1_{L_n\geq nq}
\enqs
where $L_n$ is sampled with a  default probability  $p(\theta)$ $=$ $\P_\theta[Y_k=1]$ $=$ $pe^\theta/(1-p+pe^\theta)$.  
Moreover,  see Remark \ref{remIS}, this estimator is asymptotically optimal, as $n$ goes to infinity,  for the choice of parameter $\theta_q$ $\geq$ $0$  attaining the argmax in $\theta q - \Gamma(\theta)$. 
}
\end{Remark}

\subsubsection{Dependent obligors}

We consider the case where $\rho$ $>$ $0$. Then, conditionally on the factor $Z$, the default indicators $Y_k$ are i.i.d. with Bernoulli distribution of parameter~: 
\beq 
p(Z) &=& \P[ Y_k = 1| Z] \; = \; \P[ \rho Z + \sqrt{1-\rho^2}\eps_k > -\Phi^{-1}(p) | Z ]  \nonumber \\
&=& \Phi\Big(  \frac{\rho Z + \Phi^{-1}(p)}{\sqrt{1-\rho^2}}\Big).  \label{p(z)}
\enq
Hence, by the law of large numbers, $L_n/n$ converges in law to the random variable $p(Z)$ valued in  $(0,1)$. In order that $\{L_n \geq l_n\}$ 
becomes a rare event (without being impossible) as $n$ increases, we therefore let  $l_n/n$ approach $1$ from below. We then set
\beq \label{ln}
l_n &=& n q_n, \;\;\; \mbox{ with } \; q_n < 1, \; q_n \nearrow 1 \;\; \mbox{ as } \; n \rightarrow \infty. 
\enq  
Actually, we assume that the rate of increase of $q_n$ to $1$ is of order $n^{-a}$ with $a$ $\leq$ $1$~:  
\beq \label{qn}
1 - q_n & = &  O( n^{-a}), \;\;\; \mbox{ with }  0 < a \leq 1. 
\enq

We then state the large deviations result for the large loss threshold regime.  

\begin{Theorem} \label{theolarrisk}
In the single-factor homogeneous portfolio credit risk model \reff{single}, and with large threshold $l_n$ as in \reff{ln}-\reff{qn}, we have 
\beqs
\lim_{n\rightarrow\infty} \frac{1}{\ln n}  \ln \P[L_n \geq n q_n]  & =& - a \frac{1-\rho^2}{\rho^2}.  
\enqs
\end{Theorem}

Observe that in the above theorem, we normalize by $\ln n$, indicating that the probability decays like $n^{-\gamma}$, with $\gamma$ $=$ 
$a(1-\rho^2)/\rho^2$.  We find that the decay rate is determined  by the effect of the dependence structure in the Gaussian copula model. 
When $\rho$ is small (weak dependence between sources of credit risk), large losses occur very rarely, which is formalized by a high decay rate. 
In the opposite case, this decay rate is small when $\rho$ tends to one, which means that large losses are most likely to result from systematic risk factors.

\vspace{2mm}

\noindent {\bf Proof.} 1)  We first prove the lower bound~: 
\beq \label{lowercond}
\liminf_{n\rightarrow\infty} \frac{1}{\ln n}  \ln \P[L_n \geq n q_n]  &  \geq & - a \frac{1-\rho^2}{\rho^2}.
\enq
From Bayes formula, we have
\beq 
\P[ L_n \geq nq_n ] & \geq & \P[ L_n \geq nq_n, p(Z) = q_n ] \nonumber \\
&=&    \P[ L_n \geq nq_n | p(Z) \geq  q_n ]  \;  \P[ p(Z) \geq q_n].   \label{bayes}
\enq
For any $n$ $\geq$ $1$, we define $z_n$ $\in$ $\R$ the solution to 
\beqs
p(z_n) &=& q_n, \;\;\; n \geq 1. 
\enqs
Since $p(.)$ is an increasing one to one function, we have $\{p(Z)\geq q_n\}$ $=$ $\{Z \geq z_n\}$. Moreover, observing that $L_n$ is an increasing function of $Z$, we get 
\beqs
\P[ L_n \geq nq_n | p(Z) \geq  q_n ]  &=& \P[ L_n \geq nq_n |  Z  \geq   z_n ]  \\
& \geq &   \P[ L_n \geq nq_n |  Z  =  z_n ] \; = \;  \P[ L_n \geq nq_n | p(Z) =  q_n ],  
\enqs
so that from \reff{bayes}
\beq \label{bayes2}
\P[ L_n \geq nq_n ] & \geq &  \P[ L_n \geq nq_n |  p(Z)  =  q_n ] \P[ Z \geq z_n]. 
\enq
Now given $p(Z)$ $=$ $q_n$,   $L_n$ is binomially distributed with parameters $n$ and $q_n$, and thus
\beq \label{intercond}
\P[ L_n \geq n q_n | p(Z) = q_n] & \geq & 1 - \Phi(0) = \frac{1}{2} (> 0). 
\enq
We focus on the tail probability $\P[ Z \geq z_n]$ as $n$ goes to infinity. First, observe that since $q_n$ goes to $1$, we have $z_n$ going to infinity as $n$ tends  to infinity. Furthermore,  from the expression \reff{p(z)} of $p(z)$,  the rate of decrease \reff{qn}, and using the property  that $1-\Phi(x)$ $\simeq$ $\varphi(x)/x$ as $x$ $\rightarrow$  $\infty$,  we have 
\beqs
O(n^{-a}) \; = \; 1- q_n \; = \;  1 - p(z_n) &=& 1 -   \Phi\Big(  \frac{\rho z_n + \Phi^{-1}(p)}{\sqrt{1-\rho^2}}\Big) \\
& \simeq &  \frac{\sqrt{1-\rho^2}}{\rho z_n + \Phi^{-1}(p)} \exp\Big( - \frac{1}{2} \big(\frac{\rho z_n + \Phi^{-1}(p)}{\sqrt{1-\rho^2}} \big)^2 \Big), 
\enqs
as $n$ $\rightarrow$ $\infty$, so that by taking logarithm~: 
\beqs
a \ln n -  \frac{1}{2} \frac{\rho^2 z_n^2}{1-\rho^2} -  \ln z_n &=& O(1). 
\enqs
This implies  
\beq \label{limzn}
\lim_{n\rightarrow\infty} \frac{z_n^2}{\ln n} & = & 2 a \frac{1-\rho^2}{\rho^2}. 
\enq 
By writing 
\beqs
\P[ Z \geq z_n]  &=&   \P[ z_n \leq Z \leq z_n + 1 ]  \nonumber \\
& \geq & \frac{1}{\sqrt{2\pi}} \exp\Big( - \frac{1}{2}  (z_n+1)^2   \Big),
\enqs
we deduce with \reff{limzn}
\beqs
\liminf_{n\rightarrow\infty} \frac{1}{\ln n}  \ln \P[ Z \geq z_n]  & \geq & a \frac{1-\rho^2}{\rho^2}.
\enqs
Combining with \reff{bayes2} and  \reff{intercond}, we get the required lower bound \reff{lowercond}.

\vspace{2mm}

\noindent 2) We  now focus on the upper bound
\beq \label{uppercond}
\limsup_{n\rightarrow\infty} \frac{1}{\ln n}  \ln \P[L_n \geq n q_n]  &  \leq & - a \frac{1-\rho^2}{\rho^2}.
\enq
We introduce the conditional c.g.f. of $Y_k$~: 
\beq \label{cgfz}
\Gamma(\theta,z) &=& \ln \E\big[ e^{\theta Y_k} | Z = z] \\
&=& \ln ( 1- p(z) + p(z)e^\theta). 
\enq
Then, for any $\theta$ $\geq$ $0$, we get by Chebichev's inequality, 
\beqs
\P[ L_n \geq  nq_n | Z] & \leq &  \E\big[ e^{\theta(L_n - nq_n)} | Z\big]  \; = \; e^{-n(\theta q_n - \Gamma(\theta,Z))}, 
\enqs
so that 
\beq \label{upperintercond}
\P[ L_n \geq  nq_n | Z] & \leq &  e^{-n \Gamma^*(q_n,Z)},
\enq
where 
\beqs
\Gamma^*(q,z) \; = \;  \sup_{\theta\geq 0} [ \theta q - \Gamma(\theta,z)] & =&  \left\{ 
\begin{array}{cl}
0, & \mbox{ if } \; q \leq p(z) \\
q \ln\big(\frac{q}{p(z)}\big) + (1-q) \ln\big(\frac{1-q}{1-p(z)}\big), &   \mbox{ if } \; p(z) < q \leq 1. 
\end{array}
\right. 
\enqs
By taking expectation on both sides on \reff{upperintercond}, we get
\beq \label{FnZ}
\P[ L_n \geq n q_n ] & \leq & \E\big[ e^{F_n(Z)}\big], 
\enq
where we set $F_n(z)$ $=$ $-n\Gamma^*(q_n,z)$.  Since $\rho$ $>$ $0$, the function $p(z)$ is increasing in $z$, so $\Gamma(\theta,z)$ is an increasing function of $z$ for all $\theta$ $\geq$ $0$. Hence, $F_n(z)$ is an increasing function of $z$, which is nonpositive and attains its maximum value $0$, for all $z$ s.t. $q_n$ $=$ $p(z_n)$ $\leq$ $p(z)$, i.e. $z$ $\geq$ $z_n$.     
Moreover, by differentiation, we can check that $F_n$ is a concave function of $z$.  We now introduce a change of measure. The idea is to shift the factor mean to  reduce the variance of the term $e^{F_n(Z)}$ in the r.h.s. of \reff{FnZ}. We consider the change of measure $\P_\mu$ 
that puts  the distribution  of $Z$  to $\Nc(\mu,1)$. Its likelihood ratio is given by
\beqs
\frac{d\P_\mu}{d\P} &=& \exp\big( \mu Z - \frac{1}{2}\mu^2\big),
\enqs
so that 
\beqs
\E\big[ e^{F_n(Z)}\big] &=& \E_\mu\big[ e^{F_n(Z) - \mu Z + \frac{1}{2}\mu^2}\big],
\enqs
where $\E_\mu$ denotes the expectation under $\P_\mu$. By concavity of $F_n$, we have $F_n(Z)$ $\leq$ 
$F_n(\mu) + F'_n(\mu) (Z-\mu)$, so that
\beq \label{espmu}
\E\big[ e^{F_n(Z)}\big]   &  \leq & \E_\mu\big[ e^{F_n(\mu) + (F_n'(\mu) - \mu) Z - \mu F'_n(\mu) + \frac{1}{2}\mu^2}\big]. 
\enq
We now choose $\mu$ $=$ $\mu_n$ solution to 
\beq \label{munfirst}
F_n'(\mu_n) &=& \mu_n,
\enq 
so that the term in the expectation in the r.h.s. of \reff{espmu} does not depend on $Z$, and is therefore a constant term 
(with zero-variance).  Such a $\mu_n$ exists, since, by strict concavity of the function $z$ $\rightarrow$ $F_n(z)-\frac{1}{2}z^2$,  
equation \reff{munfirst} is the first-order equation associated to the optimization problem~: 
\beqs
\mu_n & = & {\rm arg}\max_{\mu\in\R} [F_n(\mu) - \frac{1}{2}\mu^2].
\enqs
With this choice of factor mean $\mu_n$, and by inequalities \reff{FnZ}, \reff{espmu}, we get
\beq \label{Fnmun}
\P[ L_n \geq n q_n ] & \leq & e^{F_n(\mu_n) - \frac{1}{2}\mu_n^2}.
\enq
We now prove  that $\mu_n/z_n$ converges to $1$ as $n$ goes to infinity. Actually, we show that for all $\eps$ $>$ $0$, there is 
$n_0$ large enough so that for all $n$ $\geq$ $n_0$, $z_n(1-\eps)$ $<$ $\mu_n$ $<$ $z_n$.  Since 
$F_n'(\mu_n)-\mu_n$ $=$ $0$, and the function $F_n'(z)-z$ is decreasing by concavity $F_n(z)-z^2/2$, it suffices to show that 
\beq \label{derivinter}
F_n'(z_n(1-\eps)) - z_n(1-\eps) \; > \;  0 & \mbox{ and } & F_n'(z_n) - z_n \; < \;  0. 
\enq
We have 
\beqs
F_n'(z) &=& n\Big( \frac{p(z_n)}{p(z)} - \frac{1-p(z_n)}{1-p(z)} \Big) \; \varphi\Big( \frac{\rho z + \Phi^{-1}(p)}{\sqrt{1-\rho^2}} \Big) \; 
\frac{\rho}{\sqrt{1-\rho^2}}. 
\enqs
The second inequality in \reff{derivinter} holds since $F_n'(z_n)$ $=$ $0$ and $z_n$ $>$ $0$ for $q_n$ $>$ $p$, hence for $n$ large enough.  Actually, $z_n$ goes to infinity as $n$ goes to infinity from \reff{limzn}.  For the first inequality in \reff{derivinter}, we use the property that $1-\Phi(x)$ $\simeq$ $\varphi(x)/x$ as $x$ $\rightarrow$ $\infty$, so that 
\beqs
\lim_{n\rightarrow\infty} \frac{p(z_n)}{p(z_n(1-\eps))} \; = \; 1, & \mbox{ and } & 
\lim_{n\rightarrow\infty}  \frac{1-p(z_n)}{1-p(z_n(1-\eps))} \; = \; 0. 
\enqs
From \reff{limzn}, we have 
\beqs
\varphi\Big( \frac{\rho z_n(1-\eps) + \Phi^{-1}(p)}{\sqrt{1-\rho^2}} \Big) &=& 0( n^{-a(1-\eps)^2}),
\enqs
and therefore
\beqs
F_n'(z_n(1-\eps)) &=& 0(n^{1-a(1-\eps)^2}).  
\enqs
Moreover, from \reff{limzn} and as $a$ $\leq$ $1$, we have 
\beqs
z_n(1-\eps) &=& 0(\sqrt{\ln n}) \; = \;  o(n^{1-a(1-\eps)^2})
\enqs
We deduce that for $n$ large enough $F_n'(z_n(1-\eps)) - z_n(1-\eps)$ $>$ $0$ and so \reff{derivinter}. 

Finally, recalling that $F_n$ is nonpositive, and from \reff{Fnmun}, we obtain~: 
\beq \label{limmun}
\limsup_{n\rightarrow\infty} \frac{1}{\ln n}  \ln \P[L_n \geq n q_n]  &  \leq & 
- \frac{1}{2} \lim_{n\rightarrow}  \frac{\mu_n^2}{\ln n} \; = \; - \frac{1}{2} \lim_{n\rightarrow}  \frac{z_n^2}{\ln n} \; = \; - a \frac{1-\rho^2}{\rho^2}. 
\enq
\ep

\vspace{2mm}

\noindent {\bf Application~: asymptotic optimality of two-step importance sampling estimator}

\vspace{1mm}

\noindent  Consider the estimation problem of $\P[L_n \geq nq]$. We  apply a two-step importance sampling (IS) by  using  IS  conditional on the common factors $Z$ and IS to  the distribution of the factors $Z$.  Observe that conditioning on $Z$ reduces to the problem of  the independent case studied in the previous paragraph, with default probability $p(Z)$ as defined in \reff{p(z)}, and c.g.f. $\Gamma(.,Z)$ in \reff{cgfz}.  Choose $\theta_{q_n}(Z)$ $\geq$ $0$ attaining the argmax in $\theta q_n - \Gamma(\theta,Z)$, and return the estimator
\beqs
\exp(-\theta_{q_n}(Z) L_n + n\Gamma(\theta_{q_n}(Z),Z)) 1_{L_n\geq nq_n}, 
\enqs 
where $L_n$ is sampled with a default probability  $p(\theta_q(Z),Z)$ $=$ 
$p(Z)e^{\theta_{q_n}(Z)}/(1-p(Z)+p(Z)e^{\theta_{q_n}(Z)})$. This provides  an unbiased conditional estimator of $\P[L_n\geq nq_n|Z]$ and an asymptotically optimal conditional variance. We further apply IS to the factor $Z$ $\sim$ $\Nc(0,1)$ under $\P$, by shifting the factor mean to $\mu$, and then considering 
the estimator 
\beq \label{esti2IS}
\exp(-\mu Z+ \frac{1}{2}\mu^2) \exp(-\theta_{q_n}(Z) L_n + n\Gamma(\theta_{q_n}(Z),Z)) 1_{L_n\geq nq_n},
\enq 
where $Z$ is sampled from $\Nc(\mu,1)$.  To summarize,  the two-step IS estimator is generated as follows~:
\begin{itemize}
\item Sample $Z$ from $\Nc(\mu,1)$
\item Compute  $\theta_{q_n}(Z)$ and $p(\theta_{q_n}(Z),Z)$
\item Return the estimator \reff{esti2IS} where $L_n$ is sampled with default probability $p(\theta_{q_n}(Z),Z)$. 
\end{itemize} 
By construction, this provides an unbiaised estimator of $\P[L_n \geq nq_n]$, and the key point is  to specify the choice of $\mu$ in order to reduce the global variance or equivalently the second moment $M_n^2(\mu,q_n)$ 
of this estimator.  First, recall from Cauchy-Schwarz's inequality~: $M_n^2(\mu,q_n)$ $\geq$ $(\P[L_n\geq nq])^2$, so that the fastest possible  
rate of decay of $M_n^2(\mu,q_n)$ is twice the probability itself~: 
\beq \label{m2inf}
\liminf_{n\rightarrow\infty} \frac{1}{\ln n} \ln M_n^2(q_n,\mu) & \geq & 2 \lim_{n\rightarrow\infty} \frac{1}{\ln n}  \ln \P[L_n \geq n q_n] . 
\enq
To achieve this twice rate, we proceed as follows. Denoting by $\bar \E$ the expectation under the IS distribution, we have
\beqs
M_n^2(\mu,q_n) &=& \bar \E\Big[ \exp(-2\mu Z+  \mu^2) \exp(-2\theta_{q_n}(Z) L_n + 2n\Gamma(\theta_{q_n}(Z),Z)) 1_{L_n\geq nq_n} \Big] \\
&\leq &  \bar \E\Big[ \exp(-2\mu Z+  \mu^2) \exp(-2n \theta_{q_n}(Z) q_n  + 2n\Gamma(\theta_{q_n}(Z),Z))  \Big] \\
&=& \bar \E\Big[ \exp(-2\mu Z + \mu^2 + 2 F_n(Z)) \Big],
\enqs 
by definition of $\theta_{q_n}(Z)$ and $F_n(z)$ $=$ $-n\sup_{\theta\geq 0}[\theta q_n - \Gamma(\theta,z)]$ introduced in the proof of the upper bound in Theorem \ref{theolarrisk}.  As in \reff{espmu}, \reff{Fnmun}, by choosing $\mu$ $=$ $\mu_n$ solution to $F_n'(\mu_n)$ $=$ $\mu_n$, we then get 
\beqs
M_n^2(\mu_n,q_n) & \leq & \exp(2 F_n(\mu_n) - \mu_n^2) \; \leq \; \exp(-\mu_n^2),
\enqs
since $F_n$ is nonpositive.  From \reff{limmun},  this yields
\beqs
\limsup_{n\rightarrow\infty} \frac{1}{\ln n} \ln M_n^2(\mu_n,q_n) & \leq & -2 a \frac{1-\rho^2}{\rho^2} \; = \;  
2 \lim_{n\rightarrow\infty} \frac{1}{\ln n}  \ln \P[L_n \geq n q_n], 
\enqs
which proves together with  \reff{m2inf} that 
\beqs
\lim_{n\rightarrow\infty} \frac{1}{\ln n} \ln M_n^2(\mu_n,q_n) & = & -2 a \frac{1-\rho^2}{\rho^2} \; = \;  
2 \lim_{n\rightarrow\infty} \frac{1}{\ln n}  \ln \P[L_n \geq n q_n],
\enqs
and thus  the estimator \reff{esti2IS} for the choice $\mu$ $=$ $\mu_n$ is asymptotically optimal.  The choice of $\mu$ $=$ $z_n$ also leads to an asymptotically optimal estimator.

\begin{Remark}
{\rm We also prove by similar methods  large deviation results for the loss distribution in the limiting regime where individual loss probabilities 
decrease toward zero, see \cite{glakansha06} for the details. 
This setting is relevant to portfolios of highly-rated obligors, for which one-year default probabilities are extremely small. 
This is also relevant to measuring risk over short time horizons.  In this limiting regime, we set
\beqs
l_n \; = \; nq, \;\;\; \mbox { with } \; 0 < q < 1, \;\;\; p \; = \; p_n \; = \; O(e^{-na}),  \mbox { with } \;  a > 0. 
\enqs
Then, 
\beqs
\lim_{n\rightarrow\infty} \frac{1}{n} \ln \P[ L_n  \geq nq ] &=& - \frac{a}{\rho^2},
\enqs
and we may construct  similarly as in the case of large losses, a two-step IS asymptotically optimal estimator. 
}
\end{Remark}

\subsection{A large deviations approach to optimal long term investment}

\subsubsection{An asymptotic  outperforming benchmark criterion}

A popular approach for institutional managers is concerned  about the performance of their portfolio relative to the achievement 
of a given benchmark. This means that investors are interested in maximizing the probability that their wealth exceed  a predetermined index.  
Equivalently, this may be also formulated as the problem of minimizing the probability that the wealth of the investor falls below a specified value. 
This target problem was studied by several authors  for a goal achievement in finite time horizon, see e.g.  \cite{bro99} or \cite{folleu99}. 
Recently, and in a static framework, the paper  \cite{stu03}  considered an asymptotic version of this outperformance criterion when time horizon goes to infinity, which leads to a large deviations portfolio criterion.  
To illustrate the purpose, let us consider the following toy example.  Suppose that an investor trades a number $\alpha$ of shares in stock of price $S$, and keep it until time $T$. Her wealth at time $T$ is then $X_T^\alpha$ $=$ $\alpha S_T$. For simplicity, 
we take a Bachelier model for the stock price~: $S_t$ $=$ $\mu t + \sigma W_t$, where $W$ is a brownian motion.   We now look at the behavior of the average wealth when time horizon $T$ goes to infinity. 
By the law of large numbers,  for any $\alpha$ $\in$ $\R$, the average wealth converges a.s. to~: 
\beqs
\bar X_T^\alpha \; : = \; \frac{X_T^\alpha}{T} \; = \; \alpha \mu + \alpha \sigma \frac{W_T}{T} & \longrightarrow & \alpha \mu,  
\enqs
when $T$ goes to infinity.  
When considering positive stock price, as in the  Black-Scholes model, the relevant ergodic mean is the average of the growth rate, i.e. the logarithm of the wealth.  Fix some benchmark level $x$ $\in$ $\R$. Then, from Cramer's theorem, the probability of outperforming $x$  decays 
exponentially fast as~: 
\beqs
\P[ \bar X_T^\alpha  \geq x ] & \simeq &  e^{-I(x,\alpha) T},
\enqs
in the sense that $\lim_{T\rightarrow \infty} \frac{1}{T} \ln \P[\bar X_T^\alpha \geq x]$ $=$ $-I(x,\alpha)$, where 
\beqs
I(x,\alpha) &=& \sup_{\theta\in\R} [ \theta x - \Gamma(\theta,\alpha)] \\
\Gamma(\theta,\alpha) &=& \frac{1}{T} \ln \E[e^{\theta  X_T^\alpha}]. 
\enqs
Thus, the lower is the decay rate $I(x,\alpha)$, the more chance there is of realizing a portfolio performance above $x$.  
The asymptotic version of the outperforming benchmark criterion is then formulated  as~: 
\beq \label{varinfinity} 
\sup_{\alpha \in \R}  \lim_{T\rightarrow \infty} \frac{1}{T} \ln \P[\bar X_T^\alpha \geq x]  &=& - \inf_{\alpha \in\R} I(x,\alpha).
\enq
In this simple example, the quantities involved are all explicit~: 
\beqs
\Gamma(\theta,\alpha) &=&  \theta\alpha \mu + \frac{(\theta\alpha\sigma)^2}{2} \\
I(x,\alpha) &=& \left\{ \begin{array}{cl} 
                            \frac{1}{2} \Big( \frac{\alpha\mu - x}{\alpha\sigma} \Big)^2, &  \alpha \neq 0 \\
                             0, & \alpha = 0, \; x = 0 \\
                             \infty, &  \alpha = 0, \; x \neq 0. 
                             \end{array}
                             \right. 
\enqs
The solution to \reff{varinfinity} is then given by  $\alpha^*$ $=$ $x/\mu$, which means that the associated expected wealth 
$\E[\bar X_T^{\alpha^*}]$ is equal to the target $x$.

We now develop an asymptotic  dynamic version of the outperformance management criterion due to \cite{pha03a}. 
Such a problem corresponds to an ergodic objective of beating a given benchmark, and may be of particular interest for institutional managers with long term horizon, like mutual funds. On the other hand, stationary long term horizon problems are expected to be more tractable than finite horizon problems, and should provide some good insight for management problems with long, but finite, time horizon.

We formulate the problem in a rather abstract setting. Let $Z$ $=$ $(X,Y)$ be a process valued in $\R\times\R^d$,  controlled by $\alpha$, a control process valued in some subset $A$ of $\R^q$.  We denote by $\Ac$ the set of control processes.  
As usual, to alleviate notations, we omitted  the dependence of $Z$  $=$ $(X,Y)$ in $\alpha$ $\in$ $\Ac$.  
We shall then study the large deviations control problem~:  
\beq \label{deflardevcontrol}
v(x) &=& \sup_{\alpha\in\Ac}  \limsup_{T\rightarrow \infty} \frac{1}{T} \ln \P[\bar X_T \geq x], \;\;\; x \in \R,  
\enq
where $\bar X_T$ $=$ $X_T/T$. 
The variable $X$  should typically be viewed in finance as the (logarithm) of the wealth process, $Y$  are  factors on market (stock, volatility ...), and 
$\alpha$  represents the  trading  portfolio.

\subsubsection{Duality to the large deviations  control problem}

The large deviations control problem \reff{deflardevcontrol} is a non standard stochastic control problem, where the objective 
is usually formulated as an expectation of some functional to optimize.  In particular, in a Markovian continuous-time setting, we do not know if there is a dynamic programming principle and a corresponding Hamilton-Jacobi-Bellman equation for our problem.  We shall actually adopt a duality approach based on the relation relating rate function of a LDP and cumulant generating function.  The formal derivation is the following.  Given $\alpha$ $\in$ $\Ac$, if there is a LDP for 
$\bar X_T$ $=$ $X_T/T$, its rate function $I(.,\alpha)$ should be related by the Fenchel-Legendre transform~: 
\beqs
I(x,\alpha) &=& \sup_{\theta} [ \theta x - \Gamma(\theta,\alpha)],
\enqs
to the c.g.f. 
\beq \label{dualintro}
\Gamma(\theta,\alpha) &=& \limsup_{T\rightarrow\infty} \frac{1}{T} \ln \E[ e^{\theta X_T}]. 
\enq
In this case, we would get
\beqs
v(x) \; = \;  \sup_{\alpha\in\Ac}  \limsup_{T\rightarrow \infty} \frac{1}{T} \ln \P[\bar X_T \geq x] &=& - \inf_{\alpha\in\Ac} I(x,\alpha) \\
&=&   - \inf_{\alpha\in\Ac}  \sup_{\theta} [ \theta x - \Gamma(\theta,\alpha)],
\enqs
and so, provided that one could intervert infinum and supremum in the above relation (actually, the minmax theorem does not apply since  
$\Ac$ is not necessarily compact and $\alpha$ $\rightarrow$ $\theta x - \Gamma(\theta,\alpha)$ is not convex)~: 
\beq \label{reldual}
v(x) &=& - \sup_\theta [\theta x - \Gamma(\theta)],
\enq
where 
\beq \label{dualrisksen}
\Gamma(\theta) &=&  \sup_{\alpha\in\Ac} \Gamma(\theta,\alpha) \; = \;  
\sup_{\alpha\in\Ac} \limsup_{T\rightarrow\infty} \frac{1}{T} \ln \E[ e^{\theta X_T}]. 
\enq
Problem \reff{dualrisksen} is the dual problem via \reff{reldual} to the original problem \reff{deflardevcontrol}. We shall see in the next section that 
\reff{dualrisksen} can be reformulated as a risk-sensitive ergodic control problem, 
which is more tractable than \reff{deflardevcontrol} and is studied by dynamic programming methods leading in some cases to explicit calculations.

First, we show rigorously the duality relation between the large deviations control problem and the risk-sensitive control problem and how the optimal controls to the former one  are related to the latter one.  This result may be viewed as an extension of the G\"artner-Ellis theorem with control components.

\begin{Theorem} \label{thmprimal}
Suppose that there exists $\bar\theta$ $\in$ $(0,\infty]$ such that for 
all $\theta$ $\in$ $[0,\bar\theta)$, there exists a solution 
$\hat\alpha(\theta)$ $\in$ $\Ac$ to the dual problem $\Gamma(\theta)$,  with a limit in \reff{dualintro}, i.e. 
\begin{eqnarray} \label{duallim}
\Gamma(\theta) &=& \lim_{T\rightarrow\infty} 
\frac{1}{T} 
\ln \E\left[\exp\left(\theta X_T^{\hat\alpha(\theta)} \right)\right]. 
\end{eqnarray}
Suppose also that $\Gamma(\theta)$ is continuously differentiable 
on $[0,\bar\theta)$. 
Then for all $x$ $<$ $\Gamma'(\bar\theta)$ $:=$ 
$\lim_{\lambda\nearrow \bar\theta}\Gamma'(\theta)$, 
we get
\begin{eqnarray} \label{reldualite}
v(x) &=& - \sup_{\theta \in [0,\bar\theta)} \left[ \theta x - \Gamma(\theta) \right].
\end{eqnarray}
Moreover, the sequence of controls
\begin{eqnarray*}
\alpha_t^{*,n} &=& \left\{ 
                   \begin{array}{cc} 
                   \hat\alpha_t\left(\theta\left(x+\frac{1}{n}\right)\right),
                   & \;\;\; \Gamma'(0) < x < \Gamma'(\bar\theta) \\
                   \hat\alpha_t\left(\theta\left(\Gamma'(0) + 
                    \frac{1}{n}\right)\right),
                   & \;\;\; x \leq \Gamma'(0),
                  \end{array}
                  \right. 
\end{eqnarray*}
with $\theta(x)$ $\in$ $(0,\bar\theta)$ s.t. 
$\Gamma'(\theta(x))$ $=$ $x$ $\in$ $(\Gamma'(0),\Gamma'(\bar\theta))$, 
is nearly optimal in the sense that
\begin{eqnarray*}
\lim_{n\rightarrow\infty} \limsup_{T\rightarrow\infty} \frac{1}{T} 
\ln \P\left[ \bar X_T^{\alpha^{*,n}}\geq x  \right]
&=& v(x).
\end{eqnarray*}
\end{Theorem}
{\bf Proof.}

\noindent {\it Step 1.} 
Let us  consider the Fenchel-Legendre transform of the convex function 
$\Gamma$ on $[0,\bar\theta)$~:
\begin{eqnarray} \label{deffenchel}
\Gamma^*(x) &=& 
\sup_{\theta\in [0,\bar\theta)}[\theta x - \Gamma(\theta)], 
\;\;\; x \in \R.  
\end{eqnarray}
Since $\Gamma$ is $C^1$ on $[0,\bar\theta)$, it is well-known (see e.g. Lemma 2.3.9 in \cite{demzei98}) 
that the function  $\Gamma^*$ is convex, nondecreasing and satisfies~:
\begin{eqnarray} \label{Lam*}
\Gamma^*(x) &=& 
\left\{
\begin{array}{cc}
\theta(x) x - \Gamma(\theta(x)), & \mbox{if } \;   
 \Gamma'(0)  < x <  \Gamma'(\bar\theta) \\
 0, & \mbox{if } \; x \leq \Gamma'(0),  
\end{array}
\right.
\end{eqnarray}
\begin{eqnarray} \label{resultexposed} 
\theta(x) x - \Gamma^*(x) &>& \theta(x) x'  - \Gamma^*(x'), 
\;\;\; \forall  \Gamma'(0) < x < \Gamma'(\bar\theta), 
\; \forall x' \neq x,  
\end{eqnarray} 
where $\theta(x)$ $\in$ $(0,\bar\theta)$ is s.t. 
$\Gamma'(\theta(x))$ $=$ $x$ $\in$ $(\Gamma'(0),\Gamma'(\bar\theta))$. 
Moreover, $\Gamma^*$ is continuous on $(-\infty,\Gamma'(\bar\theta))$.

\vspace{2mm}

\noindent {\it Step 2~: Upper bound.}  
For all $x$ $\in$ $\R$, $\alpha$ $\in$ $\Ac$, 
an application of Chebycheff's inequality yields~:
\begin{eqnarray*}
\P[\bar X_T \geq x ] &\leq & \exp(-\theta x T) \E[\exp(\theta X_T)], 
\;\;\; \forall \; \theta \in [0,\bar\theta),
\end{eqnarray*}
and so
\begin{eqnarray*}
\limsup_{T\rightarrow \infty} \frac{1}{T} 
\ln \P[ \bar X_T \geq x ] &\leq & -\theta x + 
\limsup_{T\rightarrow \infty} \frac{1}{T} \ln  \E[\exp(\theta X_T)], 
\;\;\; \forall \; \theta \in [0,\bar\theta). 
\end{eqnarray*}
By definitions of $\Gamma$ and  $\Gamma^*$, we deduce~:
\begin{eqnarray} \label{upper}
\sup_{\alpha\in\Ac} \limsup_{T\rightarrow \infty} \frac{1}{T} 
\ln \P[\bar X_T^\alpha \geq x ] &\leq & - \Gamma^*(x). 
\end{eqnarray}

\vspace{1mm}

\noindent {\it Step 3~: Lower bound.} 
Given $x$ $<$ $\Gamma'(\bar\theta)$, 
let us define the probability 
measure $\Q_T^n$ on $(\Omega,\Fc_T)$  via~:
\begin{eqnarray} \label{defQ}
\frac{d\Q_T^n}{d\P} &=& \exp\left[ \theta(x_n) X_T^{\alpha^{*,n}} - 
\Gamma_T(\theta(x_n),\alpha^{*,n}) \right],
\end{eqnarray}  
where $x_n$ $=$ $x+1/n$ if $x$ $>$ $\Gamma'(0)$, $x_n$ $=$ 
$\Gamma'(0)+1/n$ otherwise,  $\alpha^{*,n}$ $=$ 
$\hat\alpha(\theta(x_n))$, and  
\begin{eqnarray*} 
\Gamma_T(\theta,\alpha) &=&  
\ln \E[\exp(\theta X_T^\alpha)], \;\;\; \theta \in [0,\bar\theta), \; 
\alpha \in \Ac. 
\end{eqnarray*}
Here $n$ is large enough so that $x+1/n$ $<$ $\Gamma'(\bar\theta)$.  
We now take $\eps$ $>$ $0$ small enough so that $x$ $\leq$ $x_n-\eps$ and 
$x_n+\eps$  $<$ $\Gamma'(\bar\theta)$. 
We then have~:
\begin{eqnarray*}
\frac{1}{T} \ln \P[\bar X_T^{\alpha^{*,n}} \geq x ] &\geq& 
\frac{1}{T} \ln \P\left[ x_n - \eps < \bar X_T^{\alpha^{*,n}} <  
x_n + \eps \right] \\
&=& \frac{1}{T} \ln\left(\int \frac{d\P}{d\Q_T^n} 
1_{\left\{ x_n - \eps < \bar X_T^{\alpha^{*,n}} 
<  x_n +  \eps \right\}} d\Q_T^n \right) \\
&\geq& -\theta(x_n)\left(x_n + \eps \right) 
+ \frac{1}{T} \Gamma_T(\theta(x_n),\alpha^{*,n}) \\
& & \;\;\; + \frac{1}{T} \ln \Q_T^n \left[ x_n - \eps 
< \bar X_T^{\alpha^{*,n}} < x_n +  \eps \right], 
\end{eqnarray*} 
where we use \reff{defQ} in the last inequality. 
By definition of the dual problem, this yields~:
\begin{eqnarray}
\liminf_{T\rightarrow\infty} 
\frac{1}{T} \ln \P[\bar X_T^{\alpha^{*,n}} \geq x ] &\geq& 
-\theta(x_n)\left(x_n + \eps \right) 
+ \Gamma(\theta(x_n))  \nonumber  \\
& & \;\;\; + \liminf_{T\rightarrow\infty}
\frac{1}{T} \ln \Q_T^n \left[ x_n - \eps
< \bar X_T^{\alpha^{*,n}} < x_n +  \eps  \right] \nonumber \\
&\geq& - \Gamma^*(x_n) - \theta(x_n) \eps  \nonumber \\
& & \;\;\; + \liminf_{T\rightarrow\infty}
\frac{1}{T} \ln \Q_T^n \left[ x_n - \eps 
< \bar X_T^{\alpha^{*,n}} < x_n +  \eps   \right], 
\label{interlower1} 
\end{eqnarray}
where the second inequality follows by the definition of $\Gamma^*$
(and actually holds with equality due to \reff{Lam*}). 
We now show that~: 
\begin{eqnarray} \label{inter0}
\liminf_{T\rightarrow\infty}
\frac{1}{T} \ln Q_T^n
\left[ x_n - \eps  
< \bar X_T^{\alpha^{*,n}} < x_n +  \eps   \right] &=& 0. 
\end{eqnarray}
Denote by $\tilde\Gamma_T^n$ the c.g.f.  under
$Q_T^n$ of $X_T^{\alpha^{*,n}}$.  For all $\zeta$ $\in$ $\R$, we have by \reff{defQ}~:
\begin{eqnarray*} 
\tilde\Gamma_T^n(\zeta) &:=&  
\ln E^{\Q_T^n} [\exp(\zeta X_T^{\alpha^{*,n}})] \\
&=& \Gamma_T(\theta(x_n) + \zeta,\alpha^{*,n}) - 
\Gamma_T(\theta(x_n),\alpha^{*,n}). 
\end{eqnarray*}
Therefore, by definition of the dual problem and \reff{duallim}, 
we have for all $\zeta$ $\in$ $[-\theta(x_n),\bar\theta - \theta(x_n))$~:
\begin{eqnarray} \label{inter1} 
\limsup_{T\rightarrow\infty} \frac{1}{T} \tilde\Gamma_T^n(\zeta) &\leq& 
\Gamma(\theta(x_n)+\zeta) - \Gamma(\theta(x_n)).
\end{eqnarray} 
As in part 1) of this proof, by Chebycheff's inequality, we have for
all $\zeta$ $\in$ $[0,\bar\theta-\theta(x_n))$~: 
\begin{eqnarray*}
\limsup_{T\rightarrow \infty} \frac{1}{T} 
\ln \Q_T^n \left[\bar X_T^{\alpha^{*,n}} \geq x_n + \eps \right] 
&\leq & -\zeta( x_n +  \eps)   
+ \limsup_{T\rightarrow \infty} \frac{1}{T} \tilde\Gamma_T^n(\zeta) \\
&\leq& -\zeta\left( x_n + \eps \right) 
+ \Gamma(\zeta+\theta(x_n)) - \Gamma(\theta(x_n)), 
\end{eqnarray*}
where the second inequality follows from \reff{inter1}. We deduce
\begin{eqnarray} 
\limsup_{T\rightarrow \infty} \frac{1}{T} 
\ln \Q_T^n \left[\bar X_T^{\alpha^{*,n}} \geq x_n + \eps  \right] 
&\leq& - \sup\{ \zeta \left( x_n + \eps \right) 
- \Gamma(\zeta)~: \zeta \in [\theta(x_n),\bar\theta) \} \nonumber \\
& & \;\;\; - \Gamma(\theta(x_n)) 
+ \theta(x_n)\left( x_n + \eps \right) \nonumber \\
&\leq&  - \Gamma^*\left( x_n +  \eps \right) 
- \Gamma(\theta(x_n)) + \theta(x_n)\left( x_n +  \eps \right), \nonumber \\
&=& - \Gamma^*\left( x_n + \eps \right) + 
\Gamma^*(x_n) + \eps \theta(x_n), \label{inter2}
\end{eqnarray}
where the second inequality and the last equality follow 
from \reff{Lam*}. 
Similarly, we have for all $\zeta$ $\in$ $[-\theta(x_n),0]$~:
\begin{eqnarray*}
\limsup_{T\rightarrow \infty} \frac{1}{T} 
\ln \Q_T^n \left[\bar X_T^{\alpha^{*,n}} \leq x_n - \eps \right] 
&\leq & -\zeta \left( x_n - \eps \right)
+ \limsup_{T\rightarrow \infty} \frac{1}{T} \tilde\Gamma_T^n(\zeta) \\
&\leq& -\zeta \left( x_n - \eps \right)
 + \Gamma(\theta(x_n)+\zeta) - \Gamma(\theta(x_n)), 
\end{eqnarray*}
and so~:
\begin{eqnarray} 
\limsup_{T\rightarrow \infty} \frac{1}{T} 
\ln \Q_T^n \left[\bar X_T^{\alpha^{*,n}} \leq x_n - \eps  \right] 
&\leq& - \sup\{ \zeta \left( x_n - \eps \right) 
- \Gamma(\zeta)~: \zeta \in [0,\theta(x_n)] \} \nonumber \\
& & \;\;\; - \Gamma(\theta(x_n)) 
+ \theta(x_n)\left( x_n - \eps \right) \nonumber \\
&\leq& - \Gamma^*\left( x_n - \eps \right) 
+ \Gamma^*(\theta(x_n))  -  \eps \theta(x_n). \label{inter3} 
\end{eqnarray}
By \reff{inter2}-\reff{inter3}, we then get~:
\begin{eqnarray*}  
& & \limsup_{T\rightarrow \infty} \frac{1}{T} 
\ln \Q_T^n \left[ \left\{\bar X_T^{\alpha^{*,n}}\leq x_n-\eps \right\} 
\cup \left\{\bar X_T^{\alpha^{*,n}}\geq  x_n + \eps  \right\}\right] \\
&\leq& \max \left\{ 
\limsup_{T\rightarrow \infty} \frac{1}{T} 
\ln \Q_T^n \left[\bar X_T^{\alpha^{*,n}} \geq x_n + \eps  \right]; 
\limsup_{T\rightarrow \infty} \frac{1}{T} 
\ln \Q_T^n \left[\bar X_T^{\alpha^{*,n}} \leq x_n - \eps \right] 
\right\} \\  
&\leq& \max \left\{ 
- \Gamma^*\left( x_n + \eps \right) + 
\Gamma^*(x_n) + \eps \theta(x_n); 
- \Gamma^*\left( x_n - \eps \right) 
+ \Gamma^*(\theta(x_n))  -  \eps \theta(x_n) \right\} \\
& < & 0,
\end{eqnarray*}
where the strict inequality follows from \reff{resultexposed}.
This implies that 
$\Q_T^n[\{\bar X_T^{\alpha^{*,n}}\leq x_n - \eps\}$ 
$\cup$ $\{\bar X_T^{\alpha^{*,n}}\geq x_n +\eps\}]$
$\rightarrow$ $0$ and hence  $\Q_T^n [ x_n -\eps <\bar X_T^{\alpha^{*,n}}< x_n + \eps ]$
$\rightarrow$ $1$ as $T$ goes to infinity. In particular \reff{inter0}
is satisfied, and by sending $\eps$ to zero in \reff{interlower1}, we
get~:
\begin{eqnarray*}
\liminf_{T\rightarrow\infty} 
\frac{1}{T} \ln \P[\bar X _T^{\alpha^{*,n}} \geq x ] &\geq& 
 - \Gamma^*(x_n). 
\end{eqnarray*}
By continuity of $\Gamma^*$ on $(-\infty,\Gamma'(\bar \theta))$, we 
obtain by sending $n$ to infinity and recalling that 
$\Gamma^*(x)$ $=$ $0$ $=$ $\Gamma^*(\Gamma'(0))$ for $x$ $\leq$ 
$\Gamma'(0)$~:
\begin{eqnarray*}
\liminf_{n\rightarrow\infty} \liminf_{T\rightarrow\infty} 
\frac{1}{T} \ln \P[\bar X _T^{\alpha^{*,n}} \geq x ] &\geq& 
- \Gamma^*(x).  
\end{eqnarray*}
This last inequality combined with \reff{upper} ends the proof.  
\ep

\vspace{2mm}

\begin{Remark}
{\rm Notice that in Theorem \ref{thmprimal}, the duality relation 
\reff{reldualite} holds for $x$ $<$ $\Gamma'(\bar\theta)$. When 
$\Gamma'(\bar\theta)$ $=$ $\infty$, we say that fonction 
$\Gamma$ is steep, so that \reff{reldualite} holds for all $x$ $\in$ $\R$.  
We illustrate in the next section different cases where $\Gamma$ is steep 
or not.  
}
\end{Remark}


\begin{Remark}\label{remdegree}
{\rm  In financial applications, $X_t$ is the logarithm of an investor's 
wealth $V_t^\alpha$ at time $t$, $\alpha_t$ is the proportion of wealth 
invested in $q$ risky assets $S$ and $Y$ is some economic factor 
influencing the dynamics of $S$ and the savings account $S^0$. Hence, in a diffusion model, we have 
\begin{eqnarray*}
dX_t &=&  
\left[ r(Y_t) + \alpha_t'(\mu(Y_t) - r(Y_t)e_q) - \frac{1}{2}|\alpha_t'\vartheta(Y_t)|^2 \right]dt  
+ \alpha_t'\vartheta(Y_t) dW_t,
\end{eqnarray*}
where $\mu(y)$ (resp. $\vartheta(y)$) is the rate of return (resp. volatility) 
of the risky assets, $r(y)$ is the interest rate, and $e_q$ is the unit 
vector in $\R^q$.  

Notice that the value function of the dual problem  can
be written as~:
\begin{eqnarray*}
\Gamma(\theta) &=& \lim_{T\rightarrow \infty} \frac{1}{T} \ln
E\left[ U_\theta \left(V_T^{\hat\alpha(\theta)}\right) \right],
\end{eqnarray*}
where $U_\theta(c)$ $=$ $c^{\theta}$ is a power utility function
with Constant Relative Risk Aversion (CRRA) $1-\theta$ $>$ $0$ provided that $\theta$ $<$ $1$.  
Then, Theorem \ref{thmprimal} means that for any target level
$x$, the optimal overperformance probability of growth rate is
(approximately) directly related, for large $T$, to the
expected CRRA utility of wealth, by~:
\begin{eqnarray} \label{discuss}
P[\bar X_T^{\alpha^*} \geq x ] & \approx &
E\left[ U_{\theta(x)}
\left(V_T^{\alpha^*}\right) \right] e^{-\theta(x)xT},
\end{eqnarray}
with the convention that $\theta(x)$ $=$ $0$ for $x$ $\leq$
$\Gamma'(0)$. Hence, $1-\theta(x)$ can be interpreted as a constant
degree of relative risk aversion for an investor who has an
overperformance target level $x$. Moreover, by strict convexity of
function $\Gamma^*$ in \reff{deffenchel}, 
it is clear that $\theta(x)$ is strictly 
increasing for $x$ $>$ $\Gamma'(0)$. So an investor with a higher
target level $x$ has a lower degree of relative risk aversion
$1-\theta(x)$.  In summary, Theorem \ref{thmprimal} (or relation \reff{discuss}) inversely 
relates the target level of growth rate to the
degree of relative risk aversion in expected utility theory. 
}  
\end{Remark}

\subsubsection{Explicit calculations to the dual risk-sensitive control problem}

We now show that the dual control problem \reff{dualrisksen} may be transformed via a change of probability measure into a risk-sensitive control problem.  
We consider the framework of a general  diffusion model for $Z$ $=$ $(X,Y)$~: 
\beq
dX_t &=& b(X_t,Y_t,\alpha_t) dt + \sigma(X_t,Y_t,\alpha_t) dW_t  \;\;\; \mbox{ in } \; \R \label{edsX} \\
dY_t &=& \eta(X_t,Y_t,\alpha_t) dt + \sigma(X_t,Y_t,\alpha_t) dW_t  \;\;\; \mbox{ in } \; \R^d, \label{edsY}
\enq
where $W$ is a $m$-dimensional brownian motion on a filtered probability space $(\Omega,\Fc,\F=(\Fc_t)_{t\geq 0},\P)$, and 
$\alpha$ $=$ $(\alpha_t)_{t\geq 0}$, the control process, is $\F$-adapted and valued in some subset $A$ of $\R^q$.  We denote $\Ac$ the set of control processes. The coefficients $b$, $\eta$, 
$\sigma$ and $\gamma$ are measurable functions of their arguments, and given $\alpha$ $\in$ $\Ac$ and an initial condition, 
we assume the existence and uniqueness of a strong solution to \reff{edsX}-\reff{edsY}, which we also write  by setting $Z$ $=$ $(X,Y)$~: 
\beq \label{edsZ}
dZ_t &=& B(Z_t,\alpha_t) dt + \Sigma(Z_t,\alpha_t) dW_t. 
\enq 
From the dynamics of $X$ in \reff{edsX}, 
we may rewrite the Laplace transform of $X_T$ as~:  
\begin{eqnarray}
\E\left[\exp\left(\theta X_T \right)\right] &=& 
e^{\theta X_0} \E\left[\exp\left( \theta \int_0^T b(Z_t,\alpha_t) dt + 
\theta \int_0^T \sigma(Z_t,\alpha_t) dW_t \right)\right] \nonumber \\
&=& e^{\theta X_0} \E\left[ \xi_T^\alpha(\theta) 
\exp\left(\int_0^T \ell(\theta,Z_t,\alpha_t) dt \right) \right], 
\label{transformX}
\end{eqnarray}
where 
\begin{eqnarray*}
\ell(\theta,z,a) &=& \theta b(z,a) 
+ \frac{\theta^2}{2} |\sigma(z,a)|^2,
\end{eqnarray*}
and $\xi_t^\alpha(\theta)$ is the Dol{\'e}ans-Dade exponential local 
martingale
\begin{eqnarray}
\xi_t^\alpha(\theta) &=& 
\Ec\left( \theta \int \sigma(Z_u,\alpha_u) dW_u  \right)_t \nonumber \\
&:=&   
\exp\left(\theta \int_0^t 
\sigma(Z_u,\alpha_u) dW_u  
- \frac{\theta^2}{2} \int_0^t |\sigma(Z_u,\alpha_u)|^2 du \right), 
\;\; t\geq 0. \label{doleansxi}
\end{eqnarray} 
If $\xi^\alpha(\theta)$ is a ``true'' martingale, it defines a 
probability measure $\Q$ under which, by Girsanov's theorem, the dynamics 
of $Z$ is given by~:
\begin{eqnarray*}
dZ_t &=& G(\theta,Z_t,\alpha_t) dt + \Sigma(Z_t,\alpha_t) dW^{\Q}_t,
\end{eqnarray*}  
where $W^{\Q}$ is a $\Q$-Brownian motion and 
\begin{eqnarray*}
G(\theta,z,a) &=& \left(
                   \begin{array}{c}
                    b(z,a) + \theta |\sigma(z,a)|^2 \\
                    \eta(z,a) + \theta \gamma \sigma'(z,a)
                   \end{array}
                    \right).
\end{eqnarray*}
Hence, the dual problem may be written as a stochastic control problem 
with exponential integral cost criterion~: 
\begin{eqnarray} \label{riskdual}
\Gamma(\theta) &=& \sup_{\alpha \in \Ac}\limsup_{T\rightarrow\infty} 
\frac{1}{T} 
\ln \E^{\Q}\left[\exp\left( \int_0^T \ell(\theta,Z_t,\alpha_t) dt \right)\right], 
\;\; \theta \geq 0. 
\end{eqnarray}
For fixed $\theta$, this is an ergodic risk-sensitive control problem 
which has been studied by several authors, see e.g. \cite{flemac95}, \cite{biepli04} or   \cite{ste04} in a discrete-time setting. 
It admits a dynamic programming equation~:
\begin{eqnarray} 
\Lambda(\theta) &=& \sup_{a\in A} \left[ \frac{1}{2} 
{\rm tr}\left(\Sigma\Sigma'(z,a)D^2\phi_\theta\right) + 
G(\theta,z,a).\nabla \phi_\theta  \right. \nonumber \\
& & \left. \;\;\;\;\;\;\;\;\;   + \frac{1}{2} 
\left|\Sigma'(z,a)\nabla \phi_\theta\right|^2 + \ell(\theta,z,a)
\right], \;\;  z \in \R^{d+1}. \label{pdedual}
\end{eqnarray}
The unknown is the pair $(\Lambda(\theta),\phi_\theta)$ $\in$ 
$\R\times C^2(\R^{d+1})$, and $\Lambda(\theta)$ is a candidate for 
$\Gamma(\theta)$. The above P.D.E. is formally derived by considering the finite 
horizon problem 
\begin{eqnarray*}
u_\theta(T,z) &=& \sup_{\alpha \in \Ac}
\E^{\Q}\left[\exp\left( \int_0^T \ell(\theta,Z_t,\alpha_t) dt \right)\right],
\end{eqnarray*} 
by writing the Bellman equation for this classical control problem and 
by making the logarithm transformation
\begin{eqnarray*}
\ln u_\theta(T,z) &\simeq& \Lambda(\theta) T + \phi_\theta(z), 
\end{eqnarray*}  
for large $T$. 

One can prove rigorously that a pair solution $(\Lambda(\theta),\phi_\theta)$ to the PDE \reff{pdedual} provides a 
solution $\Lambda(\theta)$ $=$ $\Gamma(\theta)$  to the dual problem \reff{dualrisksen}, with an optimal control given by the argument max in \reff{pdedual}.  This is called a verification theorem in stochastic control theory. 
Actually, there may have multiple solutions $\phi_\theta$ to \reff{pdedual} (even up to a constant), and we need some ergodicity condition to select the good one that satisfies the verification theorem.  
We refer to \cite{pha03b} for the details, and we illustrate our purpose with an example with explicit calculations.

We consider a one-factor model where the bond price $S^0$ and the stock price $S$ evolve according to~: 
\beqs
\frac{dS_t^0}{S_t^0} \; = \;  (a_0 +  b_0 Y_t) dt, & & 
\frac{dS_t}{S_t} \; = \;  (a +  b Y_t) dt  + \sigma dW_t,
\enqs
with a factor $Y$ as an Ornstein-Uhlenbeck ergodic process: 
\beqs
dY_t &=& -k Y_t dt + d B_t,
\enqs
where $a_0$, $b_0$, $a$, $b$ are constants, $k$, $\sigma$ are positive constants, and $W$, $B$ are two brownian motions, supposed non correlated for simplicity.  This includes Black-Scholes, Platen-Rebolledo or Vasicek models.  The (self-financed) wealth process $V_t$ with a proportion $\alpha_t$ invested in stock, follows the dynamics~: $dV_t$ $=$ $\alpha_t V_t \frac{dS_t}{S_t}$ $+$ $(1-\alpha_t) V_t \frac{dS_t^0}{S_t^0}$, and so the logarithm of the wealth process $X_t$ $=$ $\ln V_t$ is  governed by a linear-quadratic model~: 
\beq \label{linquaX}
dX_t &=& (\beta _0 Y_t^2 + \beta_1 \alpha_t^2 + \beta_2 Y_t\alpha_t + \beta_3 Y_t + \beta_4 \alpha_t + \beta_5) dt  
 \; + \; (\delta_ 0 Y_t + \delta_1 \alpha_t  + \delta_2) dW_t, 
\enq
where in our context, $\beta_0$ $=$ $0$, $\beta_1$ $=$ $-\sigma^2/2$, $\beta_2$ $=$ $b-b_0$, $\beta_3$ $=$ $b_0$, $\beta_4$ $=$ $a-a_0$, 
$\beta_5$ $=$ $a_0$, $\delta_0$ $=$ $0$, $\delta_1$ $=$ $\sigma$ and $\delta_2$ $=$ $0$.  Without loss of generality, we may assume that 
$\sigma$ $=$ $1$ and so $\beta_1$ $=$ $-1/2$ (embedded into $\alpha$) and $\beta_5$ $=$ $0$ (embedded into $x$).  
The P.D.E. \reff{pdedual} simplifies into  the search of a pair $(\Lambda(\theta),\phi_\theta)$ with 
$\phi_\theta$ depending only on $y$ and solution to~: 
\begin{eqnarray}
\Lambda(\theta) &=&  \frac{1}{2} \phi_\theta'' - ky \phi_\theta' + 
\frac{1}{2} |\phi_\theta'|^2  
+ \theta\left(\beta_0 + \theta \frac{\delta_0^2}{2}\right)y^2 
+ \theta(\beta_3 + \theta \delta_0\delta_2)y + \theta^2 
\frac{\delta_2^2}{2} \nonumber \\
& & \; + \frac{1}{2} \frac{\theta}{1-\theta \delta_1^2}
\left[(\beta_2+\theta\delta_0\delta_1)y + \beta_4 
+ \theta \delta_1\delta_2 \right]^2.  \label{HJBlinear}
\end{eqnarray}
Moreover, the maximum in $a$ $\in$ $\R$ of \reff{pdedual} is attained for 
\begin{eqnarray} \label{maxlinear}  
\hat \alpha(\theta,y) &=& \frac{(\beta_2+\theta\delta_0\delta_1)y + \beta_4 
+ \theta \delta_1\delta_2}{1-\theta\delta_1^2}.
\end{eqnarray}
The above calculations are valid only for $0$ $\leq$ $\theta$ $<$ $1/\delta_1^2$.  
We are looking for a quadratic solution to the ordinary differential 
equation \reff{HJBlinear}~:
\begin{eqnarray*}
\phi_\theta(y) &=& \frac{1}{2}A(\theta) y^2 + B(\theta) y.
\end{eqnarray*}
By substituting into \reff{HJBlinear}, and cancelling terms in $y^2$, $y$ and 
constant terms, we obtain
\begin{itemize}
\item a polynomial second degree equation for $A(\theta)$

\item a linear equation for $B(\theta)$, given $A(\theta)$

\item $\Lambda(\theta)$ is then expressed explicitly in function of 
$A(\theta)$ and $B(\theta)$ from \reff{HJBlinear}.
\end{itemize}
The existence of a solution to the second degree equation for $A(\theta)$, through the nonnegativity of the discriminant,  allows to determine the bound 
$\bar\theta$ and so the interval $[0,\bar\theta)$ on which $\Lambda$ is  well-defined and finite.  Moreover, we find two possible roots to the polynomial second degree equation for  $A(\theta)$, but only one satisfies the ergodicity condition.  
From Theorem \ref{thmprimal}, we deduce that 
\begin{eqnarray} \label{dualrel}  
v(x)  &=& - \sup_{\theta \in [0,\bar\theta)} \big[ \theta x - \Lambda(\theta) \big], \;\;\; 
\forall x < \Lambda'(\bar\theta), 
\end{eqnarray}
with a sequence of nearly optimal controls given by~:
\begin{eqnarray*}
\alpha_t^{*,n} &=& \left\{ 
                   \begin{array}{cc} 
                 \hat\alpha\left(\theta\left(x+\frac{1}{n}\right),Y_t\right),
                   & \;\;\; \Lambda'(0) < x < \Lambda'(\bar\theta) \\
                   \hat\alpha\left(\theta\left(\Lambda'(0) + 
                    \frac{1}{n}\right),Y_t\right),
                   & \;\;\; x \leq \Lambda'(0),
                  \end{array}
                  \right. 
\end{eqnarray*}
with $\theta(x)$ $\in$ $(0,\bar\theta)$ s.t.  $\Lambda'(\theta(x))$ $=$ $x$.  In the one-factor model described above, the function $\Lambda$ is steep, i.e.  $\Lambda'(\bar\theta)$ $=$ $\infty$, and so \reff{dualrel} holds for all $x$ $\in$ $\R$.  For example, in the Black-Scholes model, i.e. 
$b_0$ $=$ $b$ $=$ $0$,  we obtain 
\beqs
\Gamma(\theta) \; = \; \Lambda(\theta) &=&  \frac{1}{2}\frac{\theta}{1-\theta} \Big(\frac{a-a_0}{\sigma^2}\Big)^2, \;\;\; \mbox{ for } \; 
\theta \; < \bar\theta = 1, 
\enqs
\beqs
v(x) \; = \;  - \sup_{\theta\in [0,1)} [\theta x - \Gamma(\theta)] & = & \left\{ \begin{array}{cl}
- (\sqrt{x} - \sqrt{\bar x})^2, &  \; \mbox { if } \; x \geq \bar x:= \Gamma'(0) \; = \; \frac{1}{2}  \big(\frac{a-a_0}{\sigma^2}\big)^2 \\
0, &  \; \mbox { if } \; x < \bar x,
\end{array}
\right.
\enqs
$\theta(x)$ $=$ $1- \sqrt{\bar x/x}$ if $x$ $\geq$ $\bar x$,  and $0$ otherwise, and  
\beqs
\alpha_t^* &=& \left\{ \begin{array}{cl}
\sqrt{2x}, &  \; \mbox { if } \; x \geq \bar x \\
\frac{a-a_0}{\sigma^2}, & \; \mbox { if } \; x < \bar x.
\end{array}
\right.
\enqs
We observe that for an index value $x$ small enough, actually $x$ $<$ $\bar x$, the optimal investment for our large deviations criterion is equal to the optimal investment of the Merton's problem  for an investor with relative risk aversion one. When the value index is larger than $\bar x$, the optimal investment is increasing with $x$, with a degree of relative risk aversion $1-\theta(x)$  decreasing in $x$. 

In the more general linear-quadratic model \reff{linquaX}, $\Lambda$ may be steep or not depending on the parameters $\beta_i$ and $\delta_i$. 
We refer to \cite{pha03b} for the details.  Some variants and extensions of this large deviations control problem are studied in \cite{hatsek05}  and 
\cite{akigaukol05}.

\section{Conclusion}

\setcounter{equation}{0} \setcounter{Assumption}{0}
\setcounter{Theorem}{0} \setcounter{Proposition}{0}
\setcounter{Corollary}{0} \setcounter{Lemma}{0}
\setcounter{Definition}{0} \setcounter{Remark}{0}

In these notes, we developed some applications and  emphasized methods of large deviations in finance and insurance. 
These applications are  multiple, and our presentation is by no means  exhaustive.   There are numerous works dealing with large deviations techniques in the context of insurance, see e.g. \cite{dje93}, \cite{asmklu96}, or more recently \cite{kaatan05} and \cite{macsta06}.     
We also cite the paper \cite{aveboybusfri03}, which develops 
asymptotic formula for calculating implied volatility of index options.  Large deviation principle for backward stochastic differential equations 
is used by  \cite{astrai06}  in a setting motivated by credit risk management.   Other papers using large deviations in portfolio management are 
\cite{sor98} and \cite{congib00}.  
Some aspects of large deviations applied to problems in macroeconomics are studied in \cite{wil04}.

From a general viewpoint, questions related to extremal events are embedded into the extreme value theory, and we refer to the classical book 
\cite{embklumik97}  for a development of this subject, especially regarding applications in finance and insurance.

\end{document}